\documentclass[a4paper,11pt]{article}

\usepackage{amssymb,latexsym}
\usepackage{amsmath,amsfonts}
\usepackage[french]{babel}
\usepackage[latin1]{inputenc}

\newtheorem{theorem}{Th\'{e}or\`{e}me}[section]
\newtheorem{prop}{Proposition}[section]
\newtheorem{lemma}{Lemme}[section]
\newtheorem{cor}{Corollaire}[section]

\newtheorem{defi}{D\'{e}finition}[section]

\newtheorem{rem}{Remarque}[section]

\newenvironment{prooof}
{ \noindent {D\'{e}monstration: }}
{{~} \hfill  $\Box$ \par\medskip}

\newenvironment{prooofe}
{ \noindent {D\'{e}monstration: }}
{}

\newcommand{\real}{\mathbb R}

\newcommand{\nat}{\mathbb N}
\newcommand{\entier}{\mathbb Z}

\newcommand{\ra}{\rightarrow}
\newcommand{\lra}{\longrightarrow}

\newcommand{\lla}{\longleftarrow}
\newcommand{\Cinf}{\mathcal{C}^\infty}
\newcommand{\CinfR}[1]{\mathcal{C}^\infty(#1,\real)}

\newcommand{\CinfK}[1]{\mathcal{C}^\infty(#1,\real)}

\newcommand {\id} {{\mathrm{id}}}
\def \ann {{\scriptstyle{\rm ann}}}

\def \rlvt {{\scriptstyle{\rm rlvt}}}

\def \d {{\mathrm{d}}}

\def \Hom {\mathop{\hbox{\rm Hom}}\nolimits}
\def \Id {\mathop{\hbox{\rm id}}\nolimits}
\def \Che {\mathop{\hbox{\rm Che}}\nolimits}
\def \ot {\otimes}
\newcommand{\realL}{\real [[h]]}

\newcommand{\Al}{\mathcal{A}}
\newcommand{\Bl}{\mathcal{B}}
\newcommand{\Il}{\mathcal{I}}
\newcommand{\Ml}{\mathcal{M}}
\newcommand{\Cl}{\mathcal{C}}

\newcommand{\AlL}{\mathcal{A}[[h]]}
\newcommand{\BlL}{\mathcal{B}[[h]]}
\newcommand{\IlL}{\mathcal{I}[[h]]}

\newcommand{\Gs}{\mathfrak{G}}
\newcommand{\h}{\mathfrak{H}}
\newcommand{\g}{\mathfrak{g}}
\newcommand{\GsI}{\mathfrak{G}_{\mathcal{I}}}
\newcommand{\GsP}{\widetilde{\mathfrak{G}}}
\newcommand{\gs}{\mathfrak{g}}
\newcommand{\gsI}{\mathfrak{g}_{\mathcal{I}}}
\newcommand{\gsP}{\widetilde{\mathfrak{g}} }

\def \gu{\ve^\cdot\!\underline{\gsI^{\otimes \cdot}}}

\def \gextgsI{\ve^\cdot_{\g}\!{\g}\ot \underline{\gsI^{\otimes \cdot}}}

\def \gd{\ve^\cdot\!\underline{\GsI^{\otimes \cdot}}}
\def \somme{\mathop {\oplus}\limits}
\def \L{{\Lambda}}
\def \w{{\wedge}}
\def \ve{{\wedge}}
\def \tim {{\scriptstyle{\rm fois}}}

\newcommand{\Dop}{\mathbf{D}}

\newcommand{\cois}{co\"{\i}sotrope}

\newcommand{\dH}{\partial_H}
\newcommand{\dK}{\partial_K}
\makeatletter
\@addtoreset{equation}{section}

\makeatother

\begin{document}

\title{Star-repr\'{e}sentations sur des sous-vari\'{e}t\'{e}s
  co\"{\i}sotropes}

\author{{Martin Bordemann,}\cr
{\small{Laboratoire des Math\'{e}matiques et Applications}}\cr
{\small{Facult\'{e} des Sciences et Techniques, Universit\'{e} de Haute Alsace}}\cr
{\small{4, rue des Fr\`{e}res Lumi\`{e}re, 68093 Mulhouse, France}}\cr
{\small{e--mail: M.Bordemann@univ-mulhouse.fr}}\cr
{~}\cr
{Gr\'{e}gory Ginot${}^{(a)}$, Gilles Halbout${}^{(b)}$}\cr
{\small{Institut de Recherche Math\'{e}matique Avanc\'{e}e de Strasbourg}}\cr
{\small{Universit\'{e} Louis Pasteur et CNRS}}\cr
{\small{7, rue Ren\'{e} Descartes, 67084 Strasbourg, France }}\cr
{\small{${}^{(a)}$ \! e-mail:\,\texttt{ginot@math.u-strasbg.fr}}}\cr
{\small{${}^{(b)}$ \! e-mail:\,\texttt{halbout@math.u-strasbg.fr}}}\cr
{~}\cr
{ Hans-Christian Herbig${}^{(c)}$,  Stefan Waldmann${}^{(d)}$}\cr
{\small{Fakult\"at f\"ur Physik}}\cr
{\small{Albert-Ludwigs Universit\"{a}t Freiburg}}\cr
{\small{Hermann-Herder-Str. 2, D-79104 Freiburg i.~Br., R.F.A. }}\cr
{\small{${}^{(c)}$ \! e-mail:\,\texttt{herbig@majestix.physik.uni-freiburg.de}}}\cr
{\small{${}^{(d)}$ \! e-mail:\,\texttt{Stefan.Waldmann@physik.uni-freiburg.de}}}}

\markboth{Bordemann-Ginot-Halbout-Herbig-Waldmann}
{Repr\'{e}sentations des star-produits sur des sous-vari\'{e}t\'{e}s co\"{\i}sotropes}

\maketitle

\begin{abstract}

Soit $X$ une vari\'{e}t\'{e} de Poisson et $C$ une sous-vari\'{e}t\'{e}
co\"{\i}sotrope par rapport au crochet de Poisson.
Soit $\Il$ l'id\'{e}al des fonctions nulles sur $C$.
Dans ce travail nous
\'{e}tudions comment construire des star-produits $\star$ sur
$\Al=C^\infty(X)[[h]]$ pour lesquels $\Il[[h]]$ est un id\'{e}al \`{a} gauche
de mani\`{e}re \`{a} obtenir une repr\'{e}sentation de $(\Al,\star)$ sur
$\Bl[[h]]=\Al[[h]]/\Il[[h]]$, d\'{e}formant la repr\'{e}sentation naturelle
de $\Al$ sur $C^\infty(C)$. Pour cela nous montrons comment un tel
r\'{e}sultat
peut se d\'{e}duire d'une
g\'{e}n\'{e}ralisation de la conjecture de formalit\'{e} de Tamarkin
aux cocha\^{\i}nes compatibles avec le crochet de Poisson.

Nous d\'{e}montrons d'abord un th\'{e}or\`{e}me \`{a} la Hochschild-Kostant-Rosenberg
entre l'espace $\gsI$ dees champs de multivecteurs compatibles avec $C$
et  $\GsI$, un espace
naturel d'op\'{e}rateurs multidiff\'{e}rentiels compatibles avec $C$.
Ensuite nous mettons en \'{e}vidence l'existence d'une structure $G_\infty$
sur $\GsI$ et montrons enfin que les
obstructions \`{a} l'existence d'une formalit\'{e} sont
contr\^{o}\-l\'{e}es par des groupes de cohomologie. Dans le cas o\`{u}
$X=\real^n$ et
$C=\real^{n-\nu}$, nous explicitons et r\'{e}duisons ces groupes.
Nous conjecturons que ces obstructions sont nulles dans le
cas $\nu=1$, ce qui reprouverait, apr\`{e}s globalisation,
`\`{a} la Dolgushev-Fedosov' l'existence de star-repr\'{e}sentations dans
le cas des sous-vari\'{e}t\'{e}s de codimension $1$. Pour des codimensions
sup\'{e}rieures il n'est pas impossible qu'il y ait des obstructions,
li\'{e}es \`{a} certaines classes caract\'{e}ristiques des feuiletages qui
appara\^{\i}ssent dans le cas symplectique.
\end{abstract}

\begin{center}
\begin{minipage}{110mm}\footnotesize{\bf Keywords~:}
Deformation quantization,
star-product, homology \end{minipage}
\end{center}

\begin{center}
\begin{minipage}{110mm}\footnotesize{\bf AMS Classification~:}
Primary 16E40, 53D55, Secondary 18D50, 16S80
\end{minipage}
\end{center}

\section{Introduction}


Soit $X$ une vari\'{e}t\'{e} diff\'{e}rentiable de dimension $n$ et
$i:C\ra X$ une
sous-vari\'{e}t\'{e} ferm\'{e}e de codimension $l\leq n$.
Soit $\Al=\CinfR{X}$, $\Bl=\CinfR{C}$ et
\begin{equation*}
   \Il=\{f\in\Al~|~f(c)=0~\forall c\in C\}
\end{equation*}
l'id\'{e}al annulateur de $C$. En utilisant un voisinage tubulaire autour de $C$
et une partition de l'unit\'{e} on voit que la suite d'alg\`{e}bres commutatives
\begin{equation} \label{EqSuiteIAB}
   \{0\}\lra \Il \lra \Al \stackrel{i^*}{\lra} \Bl \lra \{0\}
\end{equation}
est exacte.
Soit $c\in C$ et
$$ T_cC^{\ann}=\{\beta\in T_cX^*~|~\beta(v)=0~\forall v\in
    T_cC\}.$$
Soit $P\in\Gamma(X,\Lambda^2 TX)$ une structure de Poisson.
La sous-vari\'{e}t\'{e} $C$ est dite {\em \cois ~par rapport \`{a} $P$} si
\begin{equation} \label{EqDefPoiCoi}
    P_c(\beta,\gamma)=0~~~\mathrm{quels~que~soient~}c\in C; \beta,\gamma
                              \in T_cC^{\ann}.
\end{equation}
Une condition alg\'{e}brique \'{e}quivalente est
\begin{equation*}
  \Il\hbox{ est une sous-alg\`{e}bre de Poisson de }\Al.
\end{equation*}
De mani\`{e}re g\'{e}n\'{e}rale, on d\'{e}finit
\begin{defi}
\label{compavect}
Un champ de multivecteurs
$P\in\Gamma(X,\Lambda^k TX)$ est dit {\em compatible avec $C$} (ou
{\em adapt\'{e} \`{a} $C$} si
\begin{equation}
    P_c(\beta_1,\ldots,\beta_k)=0~~~\mathrm{quels~que~soient~}c\in C;
              \beta_1,\ldots,\beta_k
                              \in T_cC^{\ann}.
\end{equation}
\end{defi}
Lorsque $P\in\Al=\Gamma(X,\Lambda^0 TX)$, cette condition \'{e}quivaut
\`{a} $P \in \Il$.

\medskip

Rappelons qu'un star-produit $\star$
sur la vari\'{e}t\'{e} de Poisson $(X,P)$ est
une multiplication associative $\real [[h]]$-bilin\'{e}aire sur le
$\realL$-module $\AlL$ telle que pour tous $f,g \in \Al$,
\begin{description}
\item[(i)]$f \star g = fg + h C_1(f,g) + \sum_{k\leq 2} h^k C_k(f,g)$,
\item[(ii)]$C_1(f,g)-C_1(g,f) = P(\d f,\d g)$,
\item[(iii)]les $ (C_r)_{r \geq 1}$ sont des op\'{e}rateurs bidiff\'{e}rentiels.
\item[(iv)]$C_r(f,1)=0=C_r(1,f)$ quel que soit $r\geq 1$.
\end{description}

\bigskip

Soit $P$ une structure de Poisson compatible
avec $C$ (au sens de la D\'{e}fini\-tion \ref{compavect}).
Le but de ce travail est de continuer le travail \cite{Bor03} (dans lequel
surtout le cas symplectique et les obstructions possibles sont \'{e}tudi\'{e}es)
par la construction des star-produits $\star$ sur $X$ pour
lesquels l'espace $\Il[[h]]$ est un id\'{e}al \`{a} gauche. De cette
fa\c{c}on, on obtiendra des repr\'{e}sentations de $(\AlL,*)$ sur l'espace
$\BlL=\AlL/\IlL$. Ceci correspond au ``coisotropic creed'' prononc\'{e}
par Jiang-Hua Lu \cite{Lu93}, c'est-\`{a}-dire,
la quantification des sous-vari\'{e}t\'{e}s \cois s par des id\'{e}aux
\`{a} gauche de
l'alg\`{e}bre $\Al$ d\'{e}form\'{e}e (voir \'{e}galement \cite{BHW00} pour cette approche
dans le cadre de la r\'{e}duction Marsden-Weinstein).
Ce probl\`{e}me est intimement li\'{e} \`{a} celui des quantification des
morphismes de Poisson, voir \cite{Bor03} pour des d\'{e}tails.

Soit $k$ un entier strictement positif et soient
$\mathcal{M}_1,\ldots,\mathcal{M}_k,\mathcal{M}$ des $\Al$-$\Al$-bimodules.
D\'{e}finissons
\begin{multline*}
   \Dop^k(\mathcal{M}_1,\ldots,\mathcal{M}_k;\mathcal{M})
     =  \{\phi:\mathcal{M}_1\otimes_{\real}\cdots\otimes_{\real}\mathcal{M}_k
            \ra\mathcal{M}~|\cr
\phi\mathrm{~est~}
      k\hbox{-}\mathrm{multidiff\acute{e}rentielle}\}
\end{multline*}
On \'{e}crira $\Dop^k(\mathcal{M}_1;\mathcal{M})$ au lieu de l'expression
ci-dessus lorsque tous les modules $\mathcal{M}_1,\ldots,\mathcal{M}_k$
sont \'{e}gaux \`{a} $\mathcal{M}_1$. Pour $k=0$ on pose $\Dop^0(~;\mathcal{M})
=\mathcal{M}$, et on convient que $\Dop^k$ s'annule lorsque $k\leq -1$.
$\Dop^1(\mathcal{M}_1;\mathcal{M})$ est consid\'{e}r\'{e} comme un
$\Al$-$\Al$-bimodule
par $(f\phi g)(\eta)=f\phi(g\eta)$ quels que soient $f,g\in\Al,\phi\in
\Dop^1(\mathcal{M}_1;\mathcal{M})$ et $\eta\in \mathcal{M}_1$.
Soit $\Gs=\oplus_{k\in\entier}\Gs^k$ l'espace des cocha\^{\i}nes
de Gerstenhaber
 o\`{u}

\begin{equation*}
 \Gs^k  =  \Dop^k(\Al;\Al).
\end{equation*}

Par analogie avec les champs de multivecteurs, d\'{e}finissons maintenant
les cocha\^{\i}nes compatibles dont le r\^{o}le sera
tr\`{e}s important pour la suite de notre \'{e}tude~:
\begin{defi}
\label{compacoch}
Le sous-espace $\GsI=\oplus_{k\in\entier}\GsI^k$ de $\Gs$ des
{\sl cocha\^{\i}nes compatibles} est d\'{e}fini comme suit~:
\begin{align*}
 \GsI^k &=
  \{\phi\in \Gs^k|\phi(f_1,\ldots,f_{k-1},g)\in\Il,\cr
&\hskip3cm  \forall f_1,\ldots,f_{k-1}\in\Al,  g\in\Il\}~(\hbox{pour }k\geq 1), \cr
 \GsI^0 & =  \Il,\cr
 \GsI^k & =  \{0\}\hbox{ pour }k\leq -1.
\end{align*}
\end{defi}

Il est alors \'{e}vident qu'un star-produit $\star$ est tel que $\IlL$ est
un id\'{e}al \`{a} gauche dans
$(\AlL,*)$ si et seulement si les op\'{e}rateurs bidiff\'{e}rentiels
$(C_r)_{r \leq 0}$ le d\'{e}finissant
appartiennent
au sous-espace $\GsI^2$. Dans ce papier, nous allons \'{e}tudier comment
la construction
de Tamarkin des morphismes de formalit\'{e} (impliquant l'existence de
star-produits) peut
se restreindre \`{a} un morphisme allant des champs
de tenseurs compatibles vers les {\sl cocha\^{\i}nes compatibles}
(et donc impliquant l'existence de star-repr\'{e}sentation).

\medskip

Pour la suite de notre travail, nous d\'{e}finissons
$\GsP=\oplus_{k\in\entier}\GsP^k$
o\`{u}
\begin{align*}
\GsP^k & = \Dop^{k-1}\big(\Al;\Dop^1(\Il;\Bl)\big)
               ~(\hbox{pour }k\geq 1), \cr
\GsP^0 & =  \Al/\Il = \Bl, \cr
\GsP^k & =  \{0\}\hbox{ pour }k\leq -1.
\end{align*}
D'autre part, on peut d\'{e}finir des analogues des sous-espaces
$\GsI$ et $\GsP$ ($\in \Gs$) dans l'espace des
champs de multivecteurs $\gs=\oplus_{k\in\entier}\gs^k$ o\`{u}
\begin{align*}
   \gs^k & =  \Gamma(X,\Lambda^k TX) ~(\hbox{pour }k\geq 0), \cr
   \gs^k & =  \{0\}\hbox{ pour }k\leq -1.
\end{align*}
L'analogue de $\GsI $ est
$\gsI=\oplus_{k\in\entier}\gsI^k$ le sous-espace de $\gs$
o\`{u}
\begin{align*}
   \gsI^k & =  \{P\in\Gamma(X,\Lambda^k TX)|
                  P_c(\beta_1,\ldots,\beta_k)=0,\cr
&\hskip3cm \forall c\in C; \beta_1,\ldots,\beta_k \in T_cC^{\ann}\}
 ~(\hbox{pour }k\geq 1), \cr
   \gsI^0 & =  \Il \cr
   \gsI^k & =  \{0\}\hbox{ pour }k\leq -1.
\end{align*}

L'analogue de $\GsP$ est l'espace $\gsP=\oplus_{k\in\entier}\gsP^k$ o\`{u}
\begin{align*}
  \gsP^k  & =  \Gamma\big(C,\Lambda^k (T_CX/TC)\big)
                    ~(\hbox{pour }k\geq 1),\cr
  \gsP^0 & =  \Bl=\Al/\Il \cr
   \gsP^k & = \{0\}\hbox{ pour }k\leq -1.
\end{align*}

\bigskip

Ce papier se pr\'{e}sente comme ceci~:
\begin{itemize}
\item Dans la Section 2, nous allons montrer un analogue du
  th\'{e}or\`{e}me
de Hochschild-Kostant-Rosenberg pour les cocha\^{\i}nes compatibles
$\GsI$, dans le cas o\`{u} $X=\real^n$. Plus pr\'{e}cis\'{e}ment,
on a les suites exactes~:
$$   \{0\}\lra \GsI \lra \Gs \lra \GsP \lra \{0\},$$
$$\hbox{et }
 \{0\}\lra \gsI \lra \gs \lra \gsP \lra \{0\}.
$$
En outre, les espaces $\GsI,\Gs$ et $\GsP$ sont des complexes dont
l'op\'{e}rateur cobord est induit par le cobord de Hochschild. La version
topologique du th\'{e}or\`{e}me de Hochschild-Kostant-Rosenberg entra\^{\i}ne que
l'alg\`{e}bre de Schouten $\gs$ est la cohomologie de $\Gs$. Nous
 montrerons que
$\gsI$ est la cohomologie de $\GsI$ et que $\gsP$ est la cohomologie de
$\GsP$.
\item Dans la Section 3, nous montrerons l'existence d'une
  structure de Gerstenhaber \`{a} homotopie pr\`{e}s
  ($G_\infty$-alg\`{e}bre) sur
$\GsI$ (nous rappellerons aussi les principales d\'{e}finitions des
structures
``\`{a} homotopie pr\`{e}s'').
\item Dans la Section 4, en utilisant le
  th\'{e}or\`{e}me
de Hochschild-Kostant-Rosenberg de la Section 2, nous montrerons
que les obstructions \`{a} la construction
de morphismes de Gerstenhaber \`{a} homotopie pr\`{e}s (morphisme
$G_\infty$) entre $\gsI$ et $\GsI$ sont donn\'{e}es par des groupes de
cohomologie de
$\left(\!{\rm Hom}({\ve^\cdot\underline{\gsI^{\otimes \cdot}}}, \gsI),\left[[-,-]_S+\w,-\right]\right)$
\item Dans la Section 5, nous montrerons que l'existence de tels
morphismes
$G_\infty$ entre $\gsI$ et $\GsI$ implique l'existence de
star-repr\'{e}sentation. Plus pr\'{e}cis\'{e}ment, si on se donne
un champ de tenseurs de Poisson compatible, on peut construire un star-produit
d\'{e}fini par des cocha\^{\i}nes compatibles et donc $\IlL$ sera
un id\'{e}al \`{a} gauche dans
$(\AlL,*)$.
\item Enfin, dans la Section 6, nous expliciterons certains groupes
de cohomologie de
$\left(\!{\rm Hom}({\ve^\cdot\underline{\gsI^{\otimes
\cdot}}},\gsI),\left[[-,-]_S+\w,-\right]\right)$ (dans le cas o\`{u}
$X=\real^n$ et $C=\real^{n-\nu}$). Nous conjecturons que
ceux-ci sont nuls dans le cas $p=n-1$, comme le laissent supposer des
calculs
en petite dimension. Ceci prouverait l'existence de
star-repr\'{e}sentations
dans le cas $X=\real^n$ et $C=\real^{n-1}$ et ensuite
dans le cas g\'{e}n\'{e}ral d'une vari\'{e}t\'{e} de codimension $1$ en
utilisant un proc\'{e}d\'{e} de globalisation `\`{a} la
 Dolgushev-Fedosov' \cite{Dol}.
\item Finalement, nous discuterons des obstructions possibles
 qui appara\^{\i}ssent
  dans le cas symplectique comme des classes caract\'{e}ristiques d'Atiyah-Molino
  des feuilletages \cite{Bor03}.
\end{itemize}

\smallskip

\noindent{\bf Notations~:}
Pour deux vari\'{e}t\'{e}s
diff\'{e}rentiables $X$ et $X'$, $\Cinf(X,X')$ d\'{e}signe l'ensemble
de toutes les applications de classe $\Cinf$ de $X$ dans $X'$. Pour un
fibr\'{e} vectoriel $E$ sur une vari\'{e}t\'{e} diff\'{e}rentiable $X$, on \'{e}crira
$\Gamma(X,E)$ pour l'espace de toutes les sections de classe $\Cinf$ du
fibr\'{e} $E$.

\smallskip

\noindent{\bf Remerciments~:}
Nous remercions EUCOR pour des aides financi\`{e}res qui ont rendu possible
ce travail entre Freiburg, Strasbourg et Mulhouse.

\section{Th\'{e}or\`{e}mes de Hochschild-Kostant-Rosenberg}



Dans cette section nous allons montrer un analogue du th\'{e}or\`{e}me
de Hochschild-Kostant-Rosenberg dans le cas o\`{u}
$X=\real^n$ et $C=\real^{n-\nu}$. Dans les deux premi\`{e}res
sous-sections
nous rappellerons les propri\'{e}t\'{e}s alg\'{e}briques de
$\Gs$, $\GsI$, et $\GsP$, puis de $\gs$, $\gsI$ et
$\gsP$. Puis dans la troisi\`{e}me nous donnerons des homotopies
explicites
entres les resolutions bar et de Koszul qui nous permettrons de
prouver
notre r\'{e}sultat principal dans la quatri\`{e}me sous-section.

\subsection{Propri\'{e}t\'{e}s alg\'{e}briques de $\Gs$, $\GsI$ et $\GsP$}

\noindent
Rappelons quelques op\'{e}rations d\'{e}finies sur l'espace des cocha\^{\i}nes de
Gerstenhaber~: pour
$k\leq 1,l$ entiers, $\phi\in \Gs^k$, $\psi\in \Gs^l$ et $1\leq i\leq
k$ on pose
\begin{multline} \label{EqDefCirci}
   (\phi\circ_i \psi)\big(f_1,\ldots,f_{k+l-1}\big)\cr
     = \phi\big(f_1,\ldots,f_{i-1},\psi(f_i,\ldots,f_{l+i-1}),f_{l+i},
       \ldots,f_{k+l-1}\big),
\end{multline}
et l'on d\'{e}finit
\begin{equation*}
 \phi\circ \psi = \sum_{i=1}^k(-1)^{(i-1)(l-1)}\phi\circ_i \psi,
\end{equation*}
(on pose $\phi\circ\psi=0$ lorsque $\phi\in\Gs^k$, $k\leq 0$),
et le {\em crochet de Gerstenhaber}~:
\begin{equation*}
  [\phi,\psi]_G=\phi\circ\psi -(-1)^{(k-1)(l-1)}\psi\circ\phi.
\end{equation*}
Gerstenhaber a montr\'{e} \cite{Ger63} que $(\Gs[-1],[~,~]_G)$ est une alg\`{e}bre de
Lie gradu\'{e}e (avec $(\Gs[-1])^k=\Gs^{k+1})$). Soit $\mu$
la multiplication point par point dans $\Al$. On d\'{e}finit alors
l'op\'{e}rateur cobord de Hochschild par~:
\begin{equation} \label{EqDefCobH}
    \mathsf{b}\phi= -[\phi,\mu]_G
\end{equation}
Enfin, on d\'{e}finit la multiplication $\cup$ par~:
\begin{equation*}
    (\phi\cup\psi)\big(f_1,\ldots,f_{k+l}\big)=\mu\big(\phi(f_1,\ldots,f_k),
                    \psi(f_{k+1},\ldots,f_{k+l})\big)
\end{equation*}
et il est clair que $(\Gs,\cup)$ est une alg\`{e}bre associative gradu\'{e}e.
\begin{prop} \label{PGsI}
 L'espace $\GsI$ a les propri\'{e}t\'{e}s suivantes:
 \begin{enumerate}
  \item $\mu\in\GsI^2$.
  \item $\GsI$ est un id\'{e}al \`{a} gauche de $(\Gs,\cup)$.
  \item $\GsI[-1]$
  et une sous-alg\`{e}bre de Lie gradu\'{e}e de $(\Gs[-1],[~,~]_G)$.
  \item $(\GsI,\mathsf{b})$ est un sous-complex de $(\Gs,\mathsf{b})$.
 \end{enumerate}
\end{prop}
\begin{prooofe}
\begin{enumerate}
\item C'est \'{e}vident car $\Il$ est un id\'{e}al de $\Al$.
\item Soit $\phi\in\Gs^k$, $\psi\in \GsI^l$ et $f_{k+l}\in\Il$. On a
 $\psi(f_{k+1},\ldots,f_{k+l})\in\Il$ donc
 $(\phi\cup\psi)(f_1,\ldots,f_{k+l})\in\Il$ car $\Il$ est un id\'{e}al
 de $\Al$.
\item Soient $\phi\in\GsI^k$ et $\psi\in\GsI^l$. Regardons
 $\phi\circ_i\psi$ (si $k\leq 0$ il n'y a rien \`{a} montrer) pour $1\leq
 i\leq k$. Soit $f_{k+l-1}\in\Il$. Si $1\leq i\leq k-1$, alors
 $f_{k+l-1}$ est le dernier argument de $\phi$, donc le membre
 droit de (\ref{EqDefCirci}) appartient \`{a} $\Il$ par d\'{e}finition de $\phi$,
 et $\phi\circ_i\psi\in\GsI^{k+l-1}$. Si $i=k$, alors $f_{k+l-1}$
 est le dernier argument de $\psi$ et la valeur de $\psi$ (qui appartient
 \`{a} $\Il$ par d\'{e}finition de $\psi$) est le dernier argument de $\phi$.
 Par d\'{e}finition de $\phi$, il vient que sa valeur appartient \`{a} $\Il$,
 donc $\phi\circ_k \psi\in\GsI^{k+l-1}$.
\item C'est une cons\'{e}quence de 1., 3. et  de la
 d\'{e}finition
de $\mathsf{b}$
 (\ref{EqDefCobH}).\hfill  $\Box$
\end{enumerate}
\end{prooofe}

L'espace $\GsP$ est muni de l'op\'{e}rateur de Hochschild $\widetilde{\mathsf{b}}$
usuel~: soient $\phi\in\GsP^k$, $f_1,\ldots,f_{k}\in\Al$ et $g\in\Il$,
alors
\begin{eqnarray}\label{EqDefCobHP}
 (\widetilde{\mathsf{b}}\phi)(f_1,\ldots,f_k)(g) & = &
      f_1~\phi(f_2,\ldots,f_k)(g)  \cr
          &  & +\sum_{r=1}^{k-1}(-1)^{r}
             \phi(f_1,\ldots,f_rf_{r+1},\ldots,f_{k})(g) \cr
          &&  +(-1)^k\phi(f_1,\ldots,f_{k-1})(f_kg).
\end{eqnarray}
Pour un entier $k$ consid\'{e}rons les projections canoniques
$\Xi^k:\Gs^k\ra\GsP^k$ suivantes o\`{u} $f_1,\ldots,f_{k-1}\in\Al$ et $g\in\Il$~:
\begin{equation*}
  \big(\Xi^k (\phi)\big)\big(f_1,\ldots,f_{k-1}\big)(g) =
             i^*\big(\phi(f_1,\ldots,f_{k-1},g)\big)
\end{equation*}
pour $k\geq 1$, $\Xi^0=i^*$ et $\Xi^k=0$ quel que soit $k\leq -1$.
\begin{prop}
 On a les propri\'{e}t\'{e}s suivantes~:
 \begin{enumerate}
  \item Le diagramme suivant est une suite exacte de complexes:
  \begin{equation*}
  \{0\}\lra (\GsI,\mathsf{b}) \lra (\Gs,\mathsf{b}) \stackrel{\Xi}{\lra}
    (\GsP,\widetilde{\mathsf{b}}) \lra \{0\}.
  \end{equation*}
  En particulier, $\GsP\cong \Gs/\GsI$.
  \item $\GsP$ est un module \`{a} gauche gradu\'{e} de $(\Gs,\cup)$.
  \item $\GsP[-1]$ est un module de Lie gradu\'{e} de $(\GsI[-1],[~,~]_G)$.
 \end{enumerate}
\end{prop}
\begin{prooofe}
\begin{enumerate}
\item Il est clair que le noyau de $\Xi$ est \'{e}gal \`{a} $\GsI$ et que
$\Xi$ est surjective. Montrons que $\Xi$
est un morphisme de complexes~: soient
 $\phi\in\Gs^k$, $f_1,\dots,f_k \in \Al$, $g \in \Il$, on a
 \begin{align*}
   \big(\Xi^{k+1}(\mathsf{b}\phi)\big)\big(f_1,\ldots,f_{k}\big)(g) = &
     i^*\big((\mathsf{b}\phi)(f_1,\ldots,f_{k},g)\big) \cr
           &\hskip-3cm = i^*\big(f_1\phi(f_2,\ldots,f_k,g)\big)
  +(-1)^{k+1} i^*\big(\phi(f_1,\ldots,f_{k})g\big)\cr
           &\hskip-2.6cm +\sum_{r=1}^{k-1}(-1)^{r}
                i^*\big(\phi(f_1,\ldots,f_rf_{r+1},
                               \ldots,f_{k},g)\big) \cr
  &\hskip-2.6cm+(-1)^ki^*\big(\phi(f_1,\ldots,f_{k-1},f_kg)\big) \cr
          &\hskip-3cm =
     i^*f_1~\big(\Xi^k\phi\big)(f_2,\ldots,f_k)(g)
+ 0 \hbox{ (car }i^*g=0)\cr
          &\hskip-2.6cm +\sum_{r=1}^{k-1}(-1)^{r}
                \big(\Xi^k\phi\big)(f_1,\ldots,f_rf_{r+1},
                               \ldots,f_{k})(g)  \cr
 &\hskip-2.6cm+ (-1)^k\big(\Xi^k\phi\big)(f_1,\ldots,f_{k-1})(f_kg) \cr
          &\hskip-3cm =  (\widetilde{\mathsf{b}}\Xi^k\phi)(f_1,\ldots,f_k)(g).
 \end{align*}
\item C'est une cons\'{e}quence du fait
que $\GsI$ est un id\'{e}al \`{a} gauche de $\Gs$ (Proposition
 \ref{PGsI}, 2.).
\item Puisque $\GsI[-1]$ est une sous-alg\`{e}bre de Lie gradu\'{e}e de $\Gs[-1]$
 d'apr\`{e}s la Proposition \ref{PGsI} 3., il vient que $\Gs[-1]$ et $\GsI[-1]$
 sont des modules de Lie gradu\'{e}s de $\GsI[-1]$, donc il en est le m\^{e}me
 pour leur quotient $\GsP[-1]$.\hfill  $\Box$
\end{enumerate}
\end{prooofe}

\subsection{Propri\'{e}t\'{e}s alg\'{e}briques de
 $\gs$, $\gsI$ et $\gsP$}

\noindent
\'Etudions maintenant les espaces $\gs,\gsI$ et $\gsP$. Rappelons la
d\'{e}finition du
crochet de Schouten $[-,-]_S$ sur $\Gs$~:
soient $f,g\in\Al$, $k,l\in\nat$, $X=X_1\wedge\cdots\wedge X_k\in\gs^k$ et
$Y=Y_1\wedge\cdots\wedge Y_l\in \gs^l$
($X_1,\ldots,X_k,
Y_1,\ldots,Y_l\in\Gamma(X,TX)$), alors le crochet est d\'{e}fini par~:
\begin{align*}
 [fX,gY]_S     =&  fg\sum_{i=1}^k\sum_{j=1}^l(-1)^{i+j}[X_i,Y_j]
       \wedge
        X_1\wedge\cdots\wedge X_{i-1}\wedge X_{i+1}\wedge\cdots\wedge
        X_k \cr
      & \hskip4.5cm\wedge Y_1\wedge\cdots\wedge Y_{j-1}\wedge
         Y_{j+1}\wedge\cdots\wedge Y_l \cr
      &  + f \sum_{i=1}^k (-1)^{k-i}X_i(g)
        X_1\wedge\cdots\wedge X_{i-1}\wedge X_{i+1}\wedge\cdots\wedge
        X_k\wedge Y \cr
      &  + gX\wedge \sum_{j=1}^l (-1)^{j}Y_j(g)
        Y_1\wedge\cdots\wedge Y_{j-1}\wedge X_{j+1}\wedge\cdots\wedge
        Y_l. \cr
 \end{align*}
Il est bien connu que ce crochet ne d\'{e}pend pas de la d\'{e}composition de $X$ et de $Y$
 en produit de champs de vecteurs. En outre $(\gs[-1],[-,-]_S)$ est une alg\`{e}bre de
 Lie gradu\'{e}e~; de plus, $(\gs,\wedge)$ est une alg\`{e}bre
 associative commutative gradu\'{e}e. Enfin, pour tous
$X\in\gs^k,Y\in\gs^l,Z\in\gs^l$, on a la r\`{e}gle de d\'{e}rivation~:
 \begin{equation*}
    [X,Y\wedge Z]_S = [X,Y]_S\wedge Z + (-1)^{(k-1)l}Y\wedge [X,Z]_S
 \end{equation*}
et donc $(\gs,[-,-]_S,\wedge)$ est une alg\`{e}bre de Gerstenhaber.

\smallskip

Le produit int\'{e}rieur $i$ peut s'\'{e}tendre en une action
(toujours not\'{e}e $i$) de $(\gs,\wedge)$ sur
l'espace $\Gamma(X,\Lambda
 T^*X)$~:
 soit $x=x_1\wedge\cdots\wedge x_k\in \gs^k$,
 $f \in \Al$ et $\alpha\in\Gamma(X,\Lambda
 T^*X)$, alors
 \begin{equation*}
    i(x)\alpha =i(x_1)\cdots i(x_k)\alpha
\hbox{ et }i(f)\alpha = f \alpha.
 \end{equation*}
De m\^{e}me, la d\'{e}riv\'{e}e de Lie des champs de vecteurs s'\'{e}tend en
 une action $L$
de $(\gs[-1],[-,-]_S)$ sur
 $\Gamma(X,\Lambda T^*X)$ d\'{e}finie par
 \begin{equation} \label{EqDefLSch}
    L(x)\alpha = [i(x),\d]\alpha.
 \end{equation}
Par r\'{e}currence sur le degr\'{e} de $y$, on montre ais\'{e}ment que
 \begin{equation} \label{EqLiSch}
     [L(x),i(y)]=i([x,y]),~\forall x,y\in\gs
 \end{equation}
Enfin, puisque $\d^2=0$ on a
 \begin{equation*}
     [L(x),\d]=0\hbox{ et }    [L(x),L(y)]=L([x,y]),~\forall x,y \in \gs
 \end{equation*}

\bigskip

Consid\'{e}rons maintenant les applications  $\Psi^0=i^*:\Al\ra
\Al/\Il=\Bl$
et pour $k\geq 1$, $\Psi^k:\gs^k \lra  \gsP^k$
d\'{e}finie par
$$   x_1\wedge\cdots\wedge x_k  \mapsto
                           \big(c\mapsto
                           (x_{\stackrel{1}{c}}\mathrm{~mod~}T_cC)
                           \wedge\cdots\wedge
                           (x_{\stackrel{k}{c}}\mathrm{~mod~}T_cC)\big). $$
\begin{prop}
 La suite d'espaces vectoriels sur $\real$
 \begin{equation}
\{0\}\lra \gsI^k \lra \gs^k \stackrel{\Psi^k}{\lra} \gsP^k \lra \{0\}
 \end{equation}
 est exacte. En particulier $\gsP^k\cong \gs^k/\gsI^k$.
\end{prop}
\begin{prooof}
 Il est clair que $\Psi^k$ est bien d\'{e}finie. Le cas
 $k=0$ est une cons\'{e}quence
 de la suite exacte (\ref{EqSuiteIAB}).\\
 Montrons la surjectivit\'{e} de $\Psi^k$: soit
 $\widetilde{Y}\in\gsP^k$. Soit $E\subset
 T_CX$ un sous-fibr\'{e} tel que $T_CX=TC\oplus E$, soit $C\subset V\subset
 E$ un voisinage ouvert de la section nulle de $E$, soit $C\subset
 U\subset X$ un voisinage ouvert de $C$ dans $X$ et soit $\Phi:V\ra U$ un
 diff\'{e}omorphisme tel que $\Phi|_C$ est l'application identique (le tout
 est dit un voisinage tubulaire de $C$). $E$ est visiblement isomorphe
 \`{a} $T_CX/TC$. Pour $c\in C$ on rappelle le rel\`{e}vement vertical de $Y\in E_c$ \`{a}
 $y\in E_c$: on a $Y^{\rlvt_y}=\frac{\d}{\d t}(y+tY)|_{t=0}$,
 donc $Y^{\rlvt_y}\in T_{y}E$. Ceci induit une injection
 $(~)^{\rlvt_{\Phi(y)}}:\Lambda^k (T_cX/T_cC)\ra \Lambda^k T_{\Phi(y)}U$,
 donc un
 rel\`{e}vement vertical des sections $\widetilde{Y}\mapsto
 Y'=(\widetilde{Y})^\rlvt\in \Gamma(X,\Lambda^k TU)$ d\'{e}fini par
 ${Y'}_{\Phi(y)}=(\widetilde{Y}_c)^{\rlvt_{y}}$. Soit
 $(\psi_U,1-\psi_U)$ une partition de l'unit\'{e} subordonn\'{e}e
 au r\'{e}couvrement
 ouvert $(U,X\setminus C)$ de $X$. Alors le
 champ de vecteurs $Y$ d\'{e}fini par $Y=\psi_U Y'$ sur $U$ et
 $Y=0$ sur $X\setminus U$ est un \'{e}l\'{e}ment de $\gs^k$ tel que
 $\Psi^k Y=\widetilde{Y}$.\\
 Finalement, montrons que
$\gsI^k$ est \'{e}gal au noyau de $\Psi^k$:
 pour $c\in C$ soit $\Psi^k_c:\Lambda^kT_cX\ra \Lambda^k(T_cX/T_cC)$.
 Puisque $T_cC^{\ann}\cong (T_cX/T_cC)^*$ alors
$$    \Lambda^k (T_cC^{\ann})\cong \big(\Lambda^k (T_cX/T_cC)\big)^*.
$$
 Soit $Y\in\gs^k$. Alors $Y\in\mathrm{Ker}\Psi^k$ si
et seulement si $\Psi^k_c(X_c)=0$ quel
 que soit $c\in C$. Il s'ensuit~:
 \[
   \begin{array}{ccll}
    \Psi^k_c(Y_c)=0 &\Longleftrightarrow &
                     \xi_c\big(\Psi^k_c(Y_c)\big)=0 &~~\forall\xi_c\in
                          \Lambda^k (T_cX/T_cC)^{*} \\
                       &\Longleftrightarrow &
                     \xi_c(Y_c)=0 &~~\forall\xi_c\in
                          \Lambda^k (T_cC)^{\ann} \\
                         &\Longleftrightarrow &
                     Y_c(\xi_c)=0 &~~\forall\xi_c\in
                          \Lambda^k (T_cC)^{\ann},
  \end{array}
 \]
 et donc $\gsI^k=\mathrm{Ker}\Psi^k$.
\end{prooof}

\noindent On obtient alors une autre caract\'{e}risation de $\gsI$:
\begin{prop} \label{PgsIIl}
 Soit $k$ un entier strictement positif.
 Un \'{e}l\'{e}ment $X\in\gs^k$ appartient \`{a} $\gsI^k$ si et seulement si
 $$i(X)(\d g_1\wedge\cdots\wedge \d g_k)\in\Il\mathrm{~quels~que~soient~}
                   g_1,\cdots,g_k\in\Il
 $$
\end{prop}
\begin{prooof}
 Il est clair que $\d_c g\in T_cC^{\mathrm{ann}}$ quel que soit $c\in C$ et
 quel que soit $g\in\Il$, donc la condition ci-dessus est n\'{e}cessaire.\\
 D'autre part, soit $\beta\in T_cC^{\ann}$ et soit
 $\big(U,(x^1,\ldots,x^{n-\nu},y^1,\ldots,y^\nu)\big)$ une carte de
 sous-vari\'{e}t\'{e} de $C$ autour de $c$ (c'est-\`{a}-dire
 $U\cap C=\{u\in U~|~y^1(u)=0,\ldots,y^\nu(u)=0\}$). Alors on trouve
 $\alpha_1,\ldots,\alpha_\nu\in\real$ tels que
 $\alpha=\sum_{i=1}^\nu\alpha_i \d y^i$. Soit $(\psi_U,1-\psi_U)$ une
 partition de l'unit\'{e} subordonn\'{e}e au r\'{e}couvrement ouvert
 $(U,X\setminus \{c\})$
 de $X$. On d\'{e}finit $g= \psi_U\sum_{i=1}^\nu\alpha_i y^i$ sur $U$ et $g=0$
 sur $X\setminus U$. Visiblement $g\in\Il$ et $\d_c g=\alpha$. Alors tout
 \'{e}l\'{e}ment de $T_cC^{\ann}$ se repr\'{e}sente comme $\d_c g$ pour un
 $g\in\Il$, d'o\`{u} la suffisance de la condition.
\end{prooof}

\begin{prop}
 L'espace $\gsI$ a les propri\'{e}t\'{e}s suivantes:
 \begin{enumerate}
  \item Toute structure de Poisson $P$ sur $X$ compatible avec $C$ est dans $\gsI^2$.
  \item $(\gsI,\wedge)$ est un id\'{e}al de $(\gs,\wedge)$.
  \item $(\gsI[-1]\!,\![-,-]_S)$ est une sous-alg\`{e}bre de Lie gradu\'{e}e
       de $(\gs[-1]\!,\![-,-]_S)$.
 \end{enumerate}
\end{prop}
\begin{prooofe}
\begin{enumerate}
\item Ceci est une cons\'{e}quence directe de la d\'{e}finition
 (\ref{EqDefPoiCoi}).
\item $\Psi_c=\sum_{k=0}^n\Psi^k_c:\Lambda T_cX\ra \Lambda (T_cX/T_cC)$ est un
   homomorphisme d'alg\`{e}\-bres de Grassmann, alors $\Psi:\gs\ra\gsP$ aussi,
   donc le noyau de $\Psi$, alors $\gsI$, est un id\'{e}al par rapport \`{a} la
   multiplication ext\'{e}rieure $\wedge$.
\item Soient $x\in\gsI^k$ et $y\in\gsI^l$. Si $k=l=0$ alors
  $[x,y]_S=0\in\gsI^{-1}$. Si $k=1$ et $l=0$, alors $y\in\Il$ et
  $[x,y]_S=x(\d y)\in\Il=\gsI^{0}$ d'apr\`{e}s la proposition \ref{PgsIIl}.
  On peut supposer que $k+l\geq 2$. Soient $g_1,\ldots,g_{k+l-1}\in\Il$
  et $\gamma=\d g_1\wedge\cdots\wedge \d g_{k+l-1}$.
  Puisque $x\wedge y\in \gsI^{k+l}$, alors $i(x)i(y)\gamma=0$. De plus,
  $\d\gamma=0$. Il s'ensuit, d'apr\`{e}s l'\'{e}quation (\ref{EqLiSch}), que
  $i([x,y])\gamma=[L(x),i(y)]\gamma = i(x)\d i(y)\gamma
  +(-1)^{kl}i(y)\d i(x)\gamma$.  La $k-1$ forme $i(y)\gamma$ est une somme finie
  d'expressions de la forme~: $\pm(i(y)\gamma')\gamma''$ o\`{u} $\gamma'$
  est une $l$-forme constitu\'{e}e de $l$ \'{e}l\'{e}ments parmi
  $\d g_1,\ldots,\d g_{k+l-1}$ et $\gamma''$ est le reste. D'apr\`{e}s la
  proposition \ref{PgsIIl} la fonction $g'=\pm i(y)\gamma'$ est dans
  $\Il$,
alors
  $i(x)\d i(y)\gamma=i(x)(\d g'\wedge \gamma'')$ car $\d\gamma''=0$ et la
  $k$-forme $\d g'\wedge \gamma''$ est le produit ext\'{e}rieur
de $k$ diff\'{e}rentielles
  d'\'{e}l\'{e}ments de l'id\'{e}al, donc $i(x)
  (\d g'\wedge \gamma'')\in\Il$. Le terme
  $(-1)^{kl}i(y)\d i(x)\gamma$ appartient \'{e}galement \`{a} $\Il$ par un
  raisonnement enti\`{e}rement analogue. Alors
  $[x,y]_S\in\gsI^{k+l-1}$.\hfill  $\Box$
\end{enumerate}
\end{prooofe}
\begin{prop}
 L'espace $\gsP$ a les propri\'{e}t\'{e}s suivantes~:
 \begin{enumerate}
  \item $(\gsP,\wedge)$ est une alg\`{e}bre commutative associative gradu\'{e}e.
  \item $(\gsP,\wedge)$ est un module \`{a} gauche gradu\'{e} de $(\gs,\wedge)$.
  \item $\gsP[-1]$ est un module d'alg\`{e}bre de Lie gradu\'{e}e pour
        $(\gsI[-1],[-,-]_S)$. De plus, on a
        \[
         X.(\alpha\wedge\beta)=(X.\alpha)\wedge\beta
        +(-1)^{(k-1)l}\alpha\wedge (X.\beta)
        \]
        quels que soient $X\in\gsI^k,\alpha\in\gsP^l,\beta\in\gsP$.
 \end{enumerate}
\end{prop}
\begin{prooofe}
\begin{enumerate}
\item Ceci est \'{e}vident.
\item Puisque $\gsI$ est un id\'{e}al de $(\gs,\wedge)$ et $\gsP\cong \gs/\gsI$,
   l'\'{e}nonc\'{e} est \'{e}vident.
\item $\gs[-1]$ et $\gsI[-1]$ sont \'{e}videmment des modules de
    $(\gsI,[~,~]_S)$, et il en est de m\^{e}me pour le quotient
    $\gsP[-1]$.\hfill  $\Box$
\end{enumerate}
\end{prooofe}
\begin{rem}
 Soit $P\in\gsI^2$ une structure de Poisson compatible avec $C$. Alors
 \`{a} l'aide de l'action de $P$ l'espace $\gsP$ devient une alg\`{e}bre
 associative commutative diff\'{e}rentielle gradu\'{e}e. La cohomologie de
 $\gsP$ est dite la {\em cohomologie BRST de $C$ par rapport \`{a} $P$.}
 Cette cohomologie est munie d'une structure d'alg\`{e}bre de Poisson.
\end{rem}

\subsection{Simplification du complexe bar}
\noindent
Dans ce paragraphe $X=\real^n$.

\bigskip

Commen\c{c}ons par rappeler la r\'{e}solution `topologique'
bar de l'alg\`{e}bre $\Al=\CinfK{\real^n}$:
Soit $\Al^e=\CinfK{\real^{2n}}$
et pour tout entier positif $k$
\begin{equation*}
  CH^k= \CinfK{\real^{(k+2)n}}.
\end{equation*}
Nous noterons $(a,x_1,\ldots,x_k,b)$ (o\`{u} $a,x_1,\ldots,x_k,b\in\real^n$) pour un point de $\real^{(k+2)n}$.
L'espace $CH^k$ est un $\Al^e$-module~:
\begin{equation*}
 \Al^e\!\times CH^k\! \ra CH^k\!:(f,F)\mapsto\! \big((a,x_1,\ldots,x_k,b)
                                   \mapsto
                                   f(a,b)F(a,x_1,\ldots,x_k,b)\big)
\end{equation*}
Pour $k\geq 1$, rappelons que l'op\'{e}rateur bord de Hochschild
$\dH^k:CH^k\ra CH^{k-1}$ est d\'{e}fini par
\begin{align*}
  (\dH^k F)(a,x_1,\ldots,x_{k-1},b)  = &
                      F(a,a,x_1,\ldots,x_{k-1},b)\cr
             &
             +\sum_{r=1}^{k-1}(-1)^rF(a,x_1,\ldots,x_r,x_r,
                                            \ldots,x_{k-1},b)
                          \cr
             &
             (-1)^k F(a,x_1,\ldots,x_{k-1},b,b)
\end{align*}
Il est clair que $\dH^k$ est un morphisme de $\Al^e$-modules et
$\dH^k\dH^{k+1}=0$. L'augmentation $\epsilon:CH^0=\Al^e\ra \Al$ est d\'{e}finie
par
\begin{equation} \label{EqDefEps}
   (\epsilon F)(a)= F(a,a)~~~\forall a\in\real^n.
\end{equation}
Il est bien connu que le complexe bar
\begin{equation} \label{EqCompBar}
  \{0\}\lla \Al \stackrel{\epsilon}{\lla} CH^0 \stackrel{\dH^1}{\lla}
           CH^1 \stackrel{\dH^2}{\lla} CH^2 \stackrel{\dH^3}{\lla}\cdots
              \stackrel{\dH^k}{\lla} CH^k \stackrel{\dH^{k+1}}{\lla}
              \cdots
\end{equation}
est acyclique~: en fait, soit $h_H^{-1}:\Al\ra CH^0$ la prolongation
$\real$-lin\'{e}aire
\begin{equation*}
   (h_H^{-1}f)(a,b)=f(a,a)
\end{equation*}
et, pour $k\geq 0$, $h_H^k:CH^k\ra CH^{k+1}$ l'application $\real$-lin\'{e}aire
\begin{equation*}
   (h_H^kF)(a,x_1,\ldots,x_{k+1},b)=(-1)^{k+1}F(a,x_1,\ldots,x_k,x_{k+1}).
\end{equation*}
En \'{e}crivant $\id_H^{-1}$ pour l'application identique $\Al\ra\Al$
et $\id_H^k$ pour l'application identique $CH^k\ra CH^k$ on
montre que
\begin{eqnarray*}
       \epsilon h_H^{-1} & = & \id_H^{-1}  \\
       h_H^{-1}\epsilon + \dH^1 h_H^0 & = & \id_H^0 \\
       h_H^{k-1}\dH^k + \dH^{k+1} h_H^k & = & \id_H^k~~~\forall k\geq 1
\end{eqnarray*}
ce qui entra\^{\i}ne l'acyclicit\'{e} du complexe bar (\ref{EqCompBar}).

\bigskip

Nous allons maintenant
d\'{e}finir un autre complexe (de Koszul) acyclique pour $\Al$ en tant que $\Al^e$-module~:
soit $E=\real^n$
\begin{equation*}
   CK^k = \Al^e\otimes_{\real}\Lambda^k
   E^*~~\mathrm{quel~que~soit~}k\in\entier.
\end{equation*}
\'Evidemment, chaque $CK^k$ est un $\Al^e$-module libre. Soit
$\xi:\real^{2n}\ra E$ d\'{e}fini par
\begin{equation*}
    \xi(a,b)=a-b
\end{equation*}
Pour tout entier $k$
strictement positif, on d\'{e}finit
l'op\'{e}rateur bord de Koszul $\dK^k:CK^k\ra CK^{k-1}$
par
\begin{equation*}
    \dK^k\omega = i(\xi)\omega~~\mathrm{quel~que~soit~}\omega\in CK^k.
\end{equation*}
Autrement dit, pour $e_1,\ldots,e_k\in E$ et $a,b\in\real^n$ on a
\begin{equation*}
    (\dK^k\omega)(a,b)\big(e_2,\ldots,e_k\big)
             =\omega(a,b)\big(a-b,e_2,\ldots,e_k\big).
\end{equation*}
Il est clair que les $\dK^k$ sont des morphismes de $\Al^e$-modules
et que $\dK^k\dK^{k+1}=0$ quel que soit l'entier strictement positif $k$.
Soit $\epsilon:CK^0=\Al^e\ra\Al$ l'augmentation d\'{e}finie comme dans
(\ref{EqDefEps}). Il en r\'{e}sulte le complexe de Koszul~:
\begin{equation} \label{EqCompKos}
  \{0\}\lla \Al \stackrel{\epsilon}{\lla} CK^0 \stackrel{\dK^1}{\lla}
           CK^1 \stackrel{\dK^2}{\lla} CK^2 \stackrel{\dK^3}{\lla}\cdots
              \stackrel{\dK^k}{\lla} CK^k \stackrel{\dK^{k+1}}{\lla}
              \cdots
\end{equation}
Ce complexe est acyclique~: en effet, soit $h_K^{-1}=h_H^{-1}:\Al\ra CK^0=\Al^e=CH^0$
la prolongation $\real$-lin\'{e}aire
\begin{equation*}
   (h_K^{-1}f)(a,b)=f(a,a).
\end{equation*}
Soit $e_1,\ldots,e_n$ une base de $E$ et $e^1,\ldots,e^n$ la base duale de
$E^*$. Pour $k\geq 0$ soit $h_K^k:CK^k\ra CK^{k+1}$ l'application $\real$-lin\'{e}aire
\begin{equation*}
   (h_K^k\omega)(a,b)=-\sum_{j=1}^n e^j \wedge
                               \int_0^1dt~t^k
                               \frac{\partial\omega}{\partial b^j}(a,tb+(1-t)a).
\end{equation*}
En \'{e}crivant $\id_K^{-1}$ pour l'application identique $\Al\ra\Al$
et $\id_K^k$ pour l'application identique $CK^k\ra CK^k$ on
montre que
\begin{eqnarray}
       \epsilon h_K^{-1} & = & \id_K^{-1} \nonumber \\
       h_K^{-1}\epsilon + \dK^1 h_K^0 & = & \id_K^0 \nonumber\\
       h_K^{k-1}\dK^k + \dK^{k+1} h_K^k & = & \id_K^k~~~\forall k\geq 1
       \label{EqHHomotop}
\end{eqnarray}
ce qui entra\^{\i}ne l'acyclicit\'{e} du complexe de Koszul (\ref{EqCompKos}).

\bigskip

D\'{e}finissons enfin
les applications $F^k:CK^k\ra CH^k$ par $F^0=\id_H^0=\id_K^0$
et pour tout $\omega\in CK^k$:
\begin{equation*}
    (F^k\omega)(a,x_1,\ldots,x_k,b)=\omega(a,b)\big(x_1-a,\ldots,x_k-a\big)
\end{equation*}
quels que soient $k\in\entier$ avec $k\geq 1$ et
$a,x_1,\ldots,x_k,b\in\real^{n}$. Il est clair que les $F^k$
sont des morphismes de $\Al^e$-modules, et on montre qu'ils sont des morphismes de
complexes, c'est-\`{a}-dire, quel que soit $k\in\entier$,
\begin{equation*}
     F^k\dK^{k+1} = \dH^{k+1}F^{k+1}.
\end{equation*}
Il existe \'{e}galement des applications $G^k:CH^k\ra CK^k$ qui sont des
morphismes de $\Al^e$-modules et des morphismes de complexes~: on d\'{e}finit
$G^0=\id_H^0=\id_K^0$ et pour tout $F\in CH^k$:
\begin{eqnarray}
(G^kF)(a,b)& = &\sum_{i_1,\ldots,i_k=1}^n e^{i_1}\wedge\cdots\wedge e^{i_k}
               \int_{0}^1dt_1 \int_{0}^{t_1}dt_2 \int_{0}^{t_2}dt_3 \cdots
                \int_{0}^{t_{k-1}}dt_k \nonumber \\
           &    & ~~~~\frac{\partial^k F}{\partial x^{i_1}\cdots \partial x^{i_k}}
                 \big(a,t_1a+(1-t_1)b,\ldots,t_ka+(1-t_k)b,b\big)\nonumber
                 \\
           &    & \label{EqDefGHochKos}
\end{eqnarray}
pour tout entier $k\geq 1$. Il est \'{e}vident que chaque $G^k$ est un
homomorphisme de $\Al^e$-modules, et \`{a} l'aide d'un calcul long mais direct
on montre que quel que soit $k\in\entier$,
\begin{equation*}
    G^k\dH^{k+1} = \dK^{k+1}G^{k+1}.
\end{equation*}
On peut repr\'{e}senter les deux applications $(F^k)_{k\in\entier}$ et
$(G^k)_{k\in\entier}$ dans le diagramme commutatif suivant~:
\begin{equation*}
 \begin{array}{ccccccc}
           \cdots &
              \stackrel{\dH^k}{\lla}& CH^k                       &
              \stackrel{\dH^{k+1}}{\lla}& CH^{k+1}                           &
              \stackrel{\dH^{k+2}}{\lla} &
              \cdots \\
           \cdots &
                                    & F^k\uparrow~\downarrow G^k &
                                        & F^{k+1}\uparrow~\downarrow G^{k+1} &
                                         &
              \cdots  \\
           \cdots &
              \stackrel{\dK^k}{\lla}& CK^k                      &
              \stackrel{\dK^{k+1}}{\lla}& CK^{k+1}                           &
              \stackrel{\dK^{k+2}}{\lla} &
              \cdots
 \end{array}
\end{equation*}

\begin{lemma} \label{LThetaProj}
 Avec les notations ci-dessus on a~:
 \begin{enumerate}
  \item $G^kF^k=\id_K^k$ quel que soit l'entier positif $k$.
  \item L'op\'{e}rateur $\Theta^k=F^kG^k:CH^k\ra CH^k$ est une projection, c'est-\`{a}-dire
       $\Theta^k\Theta^k=\Theta^k$ quel que soit l'entier positif $k$.
 \end{enumerate}
\end{lemma}
\begin{prooof}
  Les deux \'{e}nonc\'{e}s sont montr\'{e}s \`{a} l'aide d'un calcul direct.
\end{prooof}

\bigskip

Nous allons maintenant montrer l'existence d'homotopies $s_H^k:CH^k\ra
CH^{k+1}$: ce sont des homomorphismes de $\Al^e$-modules tels que
\begin{equation} \label{EqVereinfHomtop}
  \id_H^k-\Theta^k = \dH^{k+1}s_H^k + s_H^{k-1}\dH^k
\end{equation}
o\`{u} $s_H^{-1}=0$ et $s_H^0= 0$. Ceci se repr\'{e}sente de la fa\c{c}on
suivante:
\begin{equation*}
 \begin{array}{ccccccccc}
                                           \cdots                   &
    \stackrel{\dH^{k-1}}{\lla}  & CH^{k-1}                           &
    \stackrel{\dH^k}{\lla}     & CH^k                               &
    \stackrel{\dH^{k+1}}{\lla} & CH^{k+1}                           &
    \stackrel{\dH^{k+2}}{\lla} &           \cdots                   \\
                                       \cdots                   &
    \searrow  \hskip0.4cm\hbox{~}
             & \downarrow &
    {\searrow} \hskip0.4cm\hbox{~}
             & \downarrow     &
    \searrow  \hskip0.4cm\hbox{~}
             & \downarrow &
    \searrow \hskip0.4cm\hbox{~}    &           \cdots
    \\
                                       \cdots                   &
    {\scriptstyle s_H^{k-2}}
             & {\scriptstyle \id_H^{k-1}-\Theta^{k-1}} &
    {\scriptstyle s_H^{k-1}}
             & {\scriptstyle \id_H^{k}-\Theta^{k}}     &
    {\scriptstyle s_H^{k}}
             & {\scriptstyle \id_H^{k+1}-\Theta^{k+1}} &
    {\scriptstyle s_H^{k+1}}          &           \cdots                   \\
                                           \cdots                   &
    \hbox{~}\hskip0.4cm\searrow
             & \downarrow &
    \hbox{~}\hskip0.4cm{\searrow}
             & \downarrow     &
    \hbox{~}\hskip0.4cm\searrow
             & \downarrow &
   \hbox{~}\hskip0.4cm \searrow          &           \cdots                   \\
                                           \cdots                   &
    \stackrel{\dH^{k-1}}{\lla}  & CH^{k-1}                           &
    \stackrel{\dH^k}{\lla}     & CH^k                               &
    \stackrel{\dH^{k+1}}{\lla} & CH^{k+1}                           &
    \stackrel{\dH^{k+2}}{\lla} &           \cdots
 \end{array}
\end{equation*}
Puisque $\id_H^0-\Theta^0=0$ et
$(\id_H^k-\Theta^k)\dH^{k+1}=\dH^{k+1}(\id_H^{k+1}-\Theta^{k+1})$ on a
$\dH^1(\id_H^1-\Theta^1)=0$. On construit $s_H^1$ ``sur les g\'{e}n\'{e}rateurs''
de la mani\`{e}re suivante~: soit $F\in CH^1$. On consid\`{e}re $\widetilde{F}$ dans
$\widetilde{C}H^1=\CinfK{\real^{5n}}$ comme $\widetilde{F}(a',a,x,b,b')=F(a',x,b')$. On
prolonge $h_H^1$ et $\id_H^1-\Theta^1$ de $CH^1$ \`{a} tout \'{e}l\'{e}ment $T$ de
$\widetilde{C}H^1$ par $(h_H^1 T)(a',a,x_1,x_2,b,b')= T(a',a,x_1,x_2,b')$
(``les variables $a',b'$ ne sont pas affect\'{e}es'') et on fait de m\^{e}me avec
$\id_H^1-\Theta^1$. Ensuite on d\'{e}finit
\begin{equation*}
    (s_H^1F)(a,x_1,x_2,b)
      =(h_H^1(id_H^1-\Theta^1)\widetilde{F})(a',a,x_1,x_2,b,b')|_{a'=a,b'=b}.
\end{equation*}
Par construction, $s_H^1$ est un homomorphisme de $\Al^e$-modules. \`A
l'aide de (\ref{EqHHomotop}) on voit que
\begin{equation*}
    \id_H^1-\Theta^1 = \dH^2 s_H^1.
\end{equation*}
On continue par r\'{e}currence~: soient $0=s_H^0,s_H^1,\ldots,s_H^k$ d\'{e}j\`{a}
construits tels que (\ref{EqVereinfHomtop}) soit satisfaite jusqu'\`{a} l'ordre $k$.
 Alors on a
\[
  \dH^{k+1}\big(\id_H^{k+1}-\Theta^{k+1}-s_H^k\dH^{k+1}\big)=0,
\]
et l'expression suivante (pour $F\in CH^{k+1}$)
\begin{multline*}
    \lefteqn{(s_H^{k+1}F)(a,x_1,\ldots,x_{k+2},b)
      =} \cr
       (h_H^{k+1}(id_H^{k+1}-\Theta^{k+1}-s_H^k\dH^{k+1})\widetilde{F})
      \big(a',a,x_1,\ldots,x_{k+2},b,b'\big))|_{a'=a,b'=b}.
      \end{multline*}
est bien d\'{e}finie.
%
%
%
%
%
%
On a alors montr\'{e} le
\begin{lemma}
On peut construire des homotopies $s_H^k:CH^k\ra CH^{k+1}$ quel que soit
 l'entier positif $k$ telles que
 \begin{enumerate}
  \item Toute $s_H^k$ est un homomorphisme de $\Al^e$-modules.
  \item Toute $s_H^k$ est construite par une `suite d'op\'{e}rations qui
  consistent en des int\'{e}grales, des d\'{e}riv\'{e}es et des \'{e}valuations'.
  \item $s_H^0=0$.
  \item $\id_H^k-\Theta^k = \dH^{k+1}s_H^k + s_H^{k-1}\dH^k$ quel que soit
        $k\in \nat$.
 \end{enumerate}
\end{lemma}

\subsection{Calcul de la cohomologie}

\noindent
On consid\`{e}re toujours le cas $X=\real^n$.
Soient
$x=(x^1,\ldots,x^{n})$ les coordonn\'{e}es canoniques de $X$ et on va
\'{e}crire $x'$ pour $(x^1,\ldots,x^{n-\nu})$ et $x''$ pour
${x''}^1=x^{n-\nu+1},\ldots,{x''}^\nu=x^n$. On note $x=(x',x'')$. Soit
\begin{equation*}
    C=\{ x\in\real^n~|~x''=0\}.
\end{equation*}
Soit $g\in\Il$, alors $g(x',0)=0$. On d\'{e}finit les fonctions $g_j$
telles que~:
\begin{multline*}
   g(x',x'')=g(x',x'')-g(x',0)=\cr
    \sum_{j=1}^l{x''}^j \int_0^1\d t\frac{\partial g}{\partial {x''}^j}(x',tx'')
     =\sum_{j=1}^l{x''}^j g_j(x',x'')
\end{multline*}
et $\Il$ est un $\Al$-module de $l$ g\'{e}n\'{e}rateurs. Il est clair que
$\Al$, $\Il$ et $\Bl=\Al/\Il$ sont des $\Al^e$-modules~: soient
$F\in\Al^e$, $f\in\Al$, $g\in\Il$ et $h\in\Bl=\CinfK{\real^{n-\nu}}$, on
pose alors
\begin{eqnarray*}
     (Ff)(x) & = & F(x,x)f(x) \\
     (Fg)(x) & = & F(x,x)g(x) \\
     (Fh)(x')& = & F\big((x',0),(x',0)\big)h(x').
\end{eqnarray*}
Dans ce qui suit, $\Cl$ d\'{e}signera le $\Al^e$-module $\Al$, $\Il$ ou $\Bl$,
 et $\Ml$ d\'{e}signera $\Cl$ ou $\Dop(\Cl,\Bl)$.
On peut munir $\Dop(\Cl,\Bl)$ d'une structure de
$\Al^e$-module~: soit $F\in\Al^e$, $\varphi\in\Dop(\Cl,\Bl)$ et $g\in\Cl$.
Alors,
\begin{itemize}
\item dans le cas $\Cl=\Al$ ou $\Cl=\Il$ l'application
$\varphi$ est de la forme
\begin{multline*}
    \varphi(g)(x') =  \sum_{r,s=0}^N\sum_{i_1,\ldots,i_r=1}^{n-\nu}
                              \sum_{j_1,\ldots,j_s=1}^{\nu}
                              \varphi_{r,s}^{i_1\cdots i_r j_1\cdots j_s}(x')
                              \cr
\frac{\partial^{r+s}g}
                              {\partial {x'}^{i_1}\cdots\partial
                              {x'}^{i_r}
                               \partial {x''}^{j_1}\cdots\partial
                              {x''}^{j_s}}(x',0)
\end{multline*}
o\`{u} $\varphi_{r,s}^{i_1\cdots i_r j_1\cdots j_s}\in\Bl$ et $g\in\Cl$.
Pour un $F\in\Al^e$ on d\'{e}finit alors
\begin{multline*}
   \big((F\varphi)(g)\big)(x')
                     \sum_{r,s=0}^N\sum_{i_1,\ldots,i_r=1}^{n-\nu}
                              \sum_{j_1,\ldots,j_s=1}^{\nu}
                              \varphi_{r,s}^{i_1\cdots i_r j_1\cdots j_s}
                              (x')\cr
                                \left.
               \frac{\partial^{r+s}\big(F((a',0),(x',x''))g(x',x'')\big)}
                              {\partial {x'}^{i_1}\cdots\partial
                              {x'}^{i_r}
                               \partial {x''}^{j_1}\cdots\partial
                              {x''}^{j_s}}(x',0)\right|_{a'=x'}.
\end{multline*}
\item Dans le cas $\Cl=\Bl$, l'application
$\varphi$ est de la forme
\begin{multline*}
    \varphi(g)(x') \sum_{r}^N\sum_{i_1,\ldots,i_r=1}^{n-\nu}
                              \varphi_{r}^{i_1\cdots i_r}(x')
               \frac{\partial^{r}g}
                              {\partial {x'}^{i_1}\cdots\partial
                              {x'}^{i_r}}(x')
\end{multline*}
o\`{u} $g,\varphi_{r}^{i_1\cdots i_r}\in\Bl$ et on d\'{e}finit
$F\varphi$ comme pr\'{e}c\'{e}demment.
\end{itemize}
Soit $\phi\in\Dop^k(\Al,\Ml)$. Pour un multi-indice $I\in\nat^n$,
$I=(p_1,\ldots,p_n)$, on d\'{e}finit l'op\'{e}rateur des d\'{e}riv\'{e}es successives
suivant (o\`{u} $|I|=p_1+\cdots+p_n$):
\[
    \partial_{x^I}=\frac{\partial^{|I|}}
                       {\partial {x^1}^{p_1}\cdots \partial{x^n}^{p_n}}.
\]
Alors $\phi$ est de la forme g\'{e}n\'{e}rale
suivante~:
\begin{itemize}
\item pour $\Ml=\Al$, $\Ml=\Il$
\begin{equation*}
  \phi(f_1,\ldots,f_k)(x)  =
                    \sum_{|I_1|,\ldots,|I_k|\leq N}
                              \phi^{I_1\cdots I_r}(x)
               \big(\partial_{x^{I_1}}f_1\big)(x)\cdots
               \big(\partial_{x^{I_k}}f_k\big)(x)
\end{equation*}
o\`{u} $\phi^{I_1\cdots I_r}$ appartient \`{a} $\Al$ ou \`{a} $\Il$.
\item Si $\Ml=\Dop(\Cl,\Bl)$,
l'op\'{e}rateur $\phi$ est de la forme
\begin{multline*}
  \big(\phi(f_1,\ldots,f_k)(g)\big)(x')= \sum_{|I_1|,\ldots,|I_k|,|J|\leq N}
                              \phi^{I_1\cdots I_r J}(x')\cr
               \big(\partial_{x^{I_1}}f_1\big)(x')\cdots
               \big(\partial_{x^{I_k}}f_k\big)(x')
               \big(\partial_{x^{J}}g\big)(x')
\end{multline*}
o\`{u} $\phi^{I_1\cdots I_r J}\in\Bl$ et $g\in \Cl$.
\end{itemize}
Dans tous ces
cas $\phi$ se prolonge de fa\c{c}on naturelle en un homomorphisme de
$\Al^e$-modules $CH^k\ra \Ml$ par ($F\in CH^k$)
\begin{itemize}
\item si $\Ml=\Al$ et $\Ml=\Il$
\begin{multline*}
 \phi(F)(x) = \sum_{|I_1|,\ldots,|I_k|\leq N}
                              \phi^{I_1\cdots I_r}(x) \cr
       \left.\big(\partial_{x_1^{I_1}}\cdots
                    \partial_{x_k^{I_k}}F(a,x_1,\ldots,x_k,b)\big)
                    \right|_{a=x_1=\cdots =x_k=b=x},
\end{multline*}
\item et si $\Ml=\Dop(\Cl,\Bl)$
\begin{multline*}
 \big(\phi(F)(g)\big)(x)  : \sum_{|I_1|,\ldots,|I_k|,|J|\leq N}
                              \phi^{I_1\cdots I_r J}(x') \cr
               \left.\partial_{x_1^{I_1}}\cdots
                    \partial_{x_k^{I_k}}\partial_{b^J}
                    \big(F(a,x_1,\ldots,x_k,b)g(b)\big)
                    \right|_{a=x_1=\cdots =x_k=b=x'}.
\end{multline*}
\end{itemize}
\begin{lemma}
Soit $\phi\in\Dop^k(\Al,\Ml)$ (avec $\Ml=\Cl$ ou $\Ml=\Dop(\Cl,\Bl)$ et
 $\Cl=\Al,\Il,\Bl$). On d\'{e}signe toujours par le m\^{e}me symbole $\phi$
son prolongement
 \`{a} $\mathrm{Hom}_{\Al^e}(CH^k,\Ml)$. Soit
 $\theta^k\phi:\Al\times\cdots\times\Al\ra\Ml$ d\'{e}fini par
\begin{itemize}
\item
$(\theta^k\phi)(f_1,\ldots,f_k)=
        \phi\big(\Theta^k(1\otimes f_1\otimes\cdots\otimes f_k\otimes
        1)\big)
$
si $\Ml=\Cl$,
\item
$   \big((\theta^k\phi)(f_1,\ldots,f_k)\big)(g)=
        \phi\big(\Theta^k(1\otimes f_1\otimes\cdots\otimes f_k\otimes
        1)\big)(g)
$ si $\Ml=\Dop(\Cl,\Bl)$.
\end{itemize}
Alors $\theta^k\phi\in\Dop^k(\Al,\Ml)$.
\end{lemma}
\begin{prooof}
Bien que $\Theta^k$ contienne des int\'{e}grales (voir le lemme \ref{LThetaProj}
 et l'\'{e}quation (\ref{EqDefGHochKos})), l'\'{e}valuation
 $a=x_1=\cdots =x_k=b=x'$ fait dispara\^{\i}tre les arguments $t_1,\ldots,t_k$
 dans les d\'{e}riv\'{e}es des fonctions $f_1,\ldots,f_k$, donc les int\'{e}grales se
 r\'{e}solvent et donnent des nombres rationnels, et des d\'{e}riv\'{e}es partielles
 restent.
\end{prooof}
De la m\^{e}me fa\c{c}on on montre la
\begin{prop}
Soit $\phi\in\Dop^k(\Al,\Ml)$ o\`{u} $\Ml=\Cl$ ou $\Ml=\Dop(\Cl,\Bl)$.
 Alors
 \begin{enumerate}
  \item $(\mathsf{b}\phi)(f_1,\ldots,f_{k+1})=\phi\big(\dH^{k+1}(1\otimes
  f_1\otimes\cdots\otimes f_{k+1}\otimes 1)\big)$ est \'{e}gal \`{a} l'op\'{e}ra\-teur
  cobord de Hochschild:
  \begin{align*}
     (\mathsf{b}\phi)(f_1,\ldots,f_{k+1}) = &
      f_1~\phi(f_2,\ldots,f_k) \cr
          &
            +\sum_{r=1}^{k}(-1)^{r}
                \phi(f_1,\ldots,f_rf_{r+1},\ldots,f_{k+1}) \cr
           &
            +(-1)^{k+1}\phi(f_1,\ldots,f_{k})f_{k+1}.
  \end{align*}
  \item Soit $S^k:\Dop^k(\Al,\Ml)\ra\Dop^{k-1}(\Al,\Ml)$ d\'{e}finie par
  \[
     (S^k\phi)(f_1,\ldots,f_{k-1})=\phi\big(s_H^{k-1}
                (1\otimes f_1\otimes\cdots\otimes f_{k-1})\big).
  \]
  Alors $S^k$ est bien d\'{e}finie et l'on a
  \[
      \id^k-\theta^k = \mathsf{b} S^{k}+S^{k+1}\mathsf{b}.
  \]
 \end{enumerate}
\end{prop}
\`A l'aide de r\'{e}sultats classiques en
alg\`{e}bre homologique, on obtient le
\begin{cor}
 La cohomologie des complexes
 $\big(\big(\Dop^k(\Al,\Ml)\big)_{k\in\entier},\mathsf{b}\big)$ est isomorphe
 \`{a} la cohomologie du sous-complexe
 $\big(\big(\theta^k\Dop^k(\Al,\Ml)\big)_{k\in\entier},\mathsf{b}\big)$.
\end{cor}

\begin{lemma}
 La suite exacte de $A^e$-modules
 \[
     \{0\}\lra \Il \lra \Al \lra \Bl \lra \{0\}
 \]
 entra\^{\i}ne la suite exacte de $A^e$-modules
 \[
     \{0\}\lla \Dop(\Il,\Bl) \lla \Dop(\Al,\Bl) \lla \Dop(\Bl,\Bl) \lla \{0\}
 \]
 et finalement, la suite exacte de complexes
 \[
     \{0\}\lla \Dop^k\big(\Al,\Dop(\Il,\Bl)\big)
     \lla \Dop^k\big(\Al,\Dop(\Al,\Bl)\big) \lla
     \Dop^k\big(\Al,\Dop(\Bl,\Bl)\big)\lla \{0\}
 \]
\end{lemma}
\begin{prooof}
 Il est clair qu'un op\'{e}rateur diff\'{e}rentiel de $\Al$ dans $\Bl$ s'annule
 sur $\Il$ si et seulement s'il ne contient que des d\'{e}riv\'{e}es partielles
 des variables $x'$, et est donc une prolongation d'un op\'{e}rateur
 diff\'{e}rentiel de $\Bl$ dans $\Bl$.
\end{prooof}

Soit $SE$ l'alg\`{e}bre sym\'{e}trique sur l'espace vectoriel $E$. Sur l'espace
vectoriel gradu\'{e}
$SE\otimes_{\real}\Lambda E=\oplus_{k\in\entier}SE\otimes_{\real}\Lambda^k E$
on a un op\'{e}rateur cobord de Koszul $\delta_K$: soit $e_1,\ldots,e_n$ une
base de $E$ et $e^1,\ldots,e^n$ la base duale. Alors
\begin{equation}
    \delta_K(L\otimes T)=\sum_{j=1}^n~i(e^j)L\otimes e_j\wedge T
\end{equation}
o\`{u} $i(e^j)$ est le produit int\'{e}rieur sym\'{e}trique qui correspond \`{a} la
d\'{e}riv\'{e}e partielle dans la direction $e^j$ si l'on interpr\`{e}te $SE$ en tant
qu'alg\`{e}bre de polyn\^{o}mes sur $E^*$.
\begin{lemma}
Le complexe
$(SE\otimes_{\real}\Lambda E,\delta_K)$ est acyclique.
\end{lemma}
Soit $\Che:\Bl\otimes_{\real}SE\ra\Dop(\Al,\Bl)$ la bijection (`de
quantification par l'ordre standard') suivante~: pour $h\in\Bl$ et
$I=(p_1,\ldots,p_n)\in\nat^n$
\begin{equation} \label{EqOrdStand}
       h \otimes (e_1)^{p_1}\cdots (e_n)^{p_n} \mapsto
     h\partial_{b^I}.
\end{equation}
Pour simplifier les notations, on posera $\hat{a}=\Che(a)$.
Soit $E_1=\real^{n-\nu}$ et $E_2=\real^\nu$. On a
\begin{equation*}
   SE\otimes \Lambda E \cong (SE_1\otimes \Lambda E_1)\otimes
                              (SE_2\otimes \Lambda E_2)
\end{equation*}
et l'on a une bijection, comme celle d\'{e}finie en (\ref{EqOrdStand}) entre
$\Bl\otimes_{\real} SE_1$ et $\Dop(\Bl,\Bl)$.

\medskip

Le th\'{e}or\`{e}me suivant simplifie le calcul de la cohomologie des
complexes $\Gs,\GsI$ et $\GsP$:
\begin{theorem}
 Les applications $\real$-lin\'{e}aires suivantes~:
 \begin{enumerate}
  \item $\chi^k:\big(\Bl\otimes_{\real}SE\otimes_{\real} \Lambda^k E,\delta_K)
          \ra \big(\theta^k\Dop^k(\Al,\Dop(\Al,\Bl)),\mathsf{b}\big)$ avec
    \begin{multline*}
       \lefteqn{\chi^k(h\otimes L\otimes T)(f_1,\ldots,f_k)(g)(x')=}
        \cr
        (-1)^{k+1}h(x') i(T)\hat{L}[b]
        \left.\big(G^k(1\otimes f_1\otimes\cdots\otimes
         f_k\otimes 1)(a,b)g(b)\big)\right|_{a=b=x'}
    \end{multline*}
    o\`{u} la notation
    $\hat{L}[b]$ indique que les d\'{e}riv\'{e}es partielles de $\hat{L}$ sont par
    rapport aux variables $b^1,\ldots,b^n$.
  \item ${\chi'}^k:\big(\Bl\otimes_{\real}SE_1\otimes_{\real}
             \oplus_{r=0}^k(\Lambda^r E_2
                 \otimes \Lambda^{k-r} E_1),(\delta_K)|_{SE_1\otimes \Lambda E_1})
          \ra \big(\theta^k\Dop^k(\Al,$ $\Dop(\Bl,\Bl)),\mathsf{b}\big)$ avec
   \begin{multline*}
       \lefteqn{\chi^k(h\otimes L\otimes T)(f_1,\ldots,f_k)(g)(x')=}
        \cr
        (-1)^{k+1}h(x') i(T)\hat{L}[b']
        \left.\big(G^k(1\otimes f_1\otimes\cdots\otimes
         f_k\otimes 1)(a,b)g(b)\big)\right|_{a=b=x'}
            \end{multline*}
 \end{enumerate}
induisent des isomorphismes de complexes.
\end{theorem}
Il en r\'{e}sulte~:
\begin{theorem} Les groupes de cohomologie des complexes suivants
se simplifient comme suit~:
 \begin{enumerate}
  \item $H\Dop^k(\Al,\Dop(\Al,\Bl))\cong \{0\}$ si $k\geq 1$.
  \item $H\Dop^k(\Al,\Dop(\Bl,\Bl))\cong \Bl\otimes \Lambda^k E_2$ quel
  que soit l'entier $k$.
  \item $H\GsP^k=H\Dop^{k-1}(\Al,\Dop(\Il,\Bl))\cong H\Dop^k(\Al,\Dop(\Bl,\Bl))
         \cong \Bl\otimes \Lambda^k E_2=\gsP^k$ quel
  que soit l'entier $k$.
 \end{enumerate}
\end{theorem}
En utilisant la suite exacte longue de cohomologie r\'{e}sultant de la suite
exacte courte des complexes $\GsI\ra\Gs\ra\GsP$ et le fait que
l'homomorphisme connectant s'annule, on obtient finalement le r\'{e}sultat
souhait\'{e}~:
\begin{theorem}
\label{theoHKR}
  Quel que soit l'entier $k$:
 \begin{enumerate}
  \item $H\Gs^k\cong \Al\otimes \Lambda^kE=\gs$.
  \item $H\GsP^k\cong \Bl \otimes \Lambda^k E_2 =\gsP^k$.
  \item $H\GsI^k\cong \gsI^k$.
 \end{enumerate}
\end{theorem}

\subsection{Les applications HKR}

\noindent Dans cette sous-section nous allons construire des
quasi-isomorphismes entre les
alg\`{e}bres de Lie diff\'{e}rentielles gradu\'{e}es $(\GsI,\mathsf{b})$ et
$(\gsI,0)$. Nous verrons que
l'application HKR usuelle correspondant \`{a} l'antisym\'{e}trisation
ne convient pas et doit \^{e}tre modifi\'{e}e.

Pour un entier positif $k$ soit $\alpha^k:\Lambda^k E\ra E^{\otimes k}$
l'aapplication d'antisym\'{e}trisa\-tion usuelle:
$v_1\wedge\cdots\wedge
  v_k\mapsto \frac{1}{k!}\sum_{\sigma\in S_k}\epsilon(\sigma)
     v_{\sigma(1)}\otimes\cdots\otimes v_{(\sigma(k))}$. Soit aussi
$P^1:SE\ra E$ la projection canonique.
En utilisant le fait que $\Gs^k\cong \Al \otimes (SE)^{\otimes k}$
(\`{a} l'aide de l'application (\ref{EqOrdStand})) nous d\'{e}finissons
les deux applications suivantes
 \begin{equation}
     \psi^k_{\scriptscriptstyle HKR}:\gs^k\ra \Gs^k: f\otimes T \mapsto
     f\otimes
     \alpha^k(T),
\end{equation}
et
\begin{equation}
   \pi^k_{\scriptscriptstyle HKR}:\Gs^k\ra \gs^k:
   f\otimes L_1\otimes\cdots\otimes L_k\mapsto f\otimes P^1(L_1)\wedge
     \cdots\wedge P^1(L_k).
\end{equation}
Nous \'{e}crirons
$\psi_{\scriptscriptstyle HKR}$ (resp. $\pi_{\scriptscriptstyle HKR}$)
pour la somme de tous
les $\psi^k_{\scriptscriptstyle HKR}$ (resp. $\pi^k_{\scriptscriptstyle
HKR}$) o\`{u} $\psi^k_{\scriptscriptstyle HKR}$ et $\pi^k_{\scriptscriptstyle
HKR}$ s'annulent pour $k\leq -1$ et $k>\dim E$.
Le th\'{e}or\`{e}me de Hochschild, Kostant et Rosenberg (HKR), \ref{theoHKR},
permet de d\'{e}duire le suivant
\begin{theorem}\label{applHKRusuel}
 Quel que soit l'entier positif $k$, les applications
 $\psi^k_{\scriptscriptstyle HKR}$ et
 $\pi^k_{\scriptscriptstyle HKR}$ ont les propri\'{e}t\'{e}s suivantes:
 \begin{enumerate}
  \item $\pi^k_{\scriptscriptstyle HKR}\psi^k_{\scriptscriptstyle HKR}=
          \mathrm{id}_{\gs^k}$.
  \item $\mathsf{b}\psi^k_{\scriptscriptstyle HKR}=0$.
  \item Soit $\phi\in\Gs^k$. Si $\mathsf{b}\phi=0$ et
      $\pi^k_{\scriptscriptstyle HKR}\phi=0$, alors $\phi$ est un cobord,
      i.e. il existe $\psi\in\Gs^{k-1}$ tel que $\phi=\mathsf{b}\psi$.
  \item Soient $X,Y\in\gs$, alors il existe $\phi\in \Gs$ tel que
      $[\psi_{\scriptscriptstyle HKR}X,\psi_{\scriptscriptstyle HKR}Y]_G
      -\psi_{\scriptscriptstyle HKR}[X,Y]_S=\mathsf{b}\phi$.
 \end{enumerate}
\end{theorem}
Malheureusement, l'application HKR $\psi_{\scriptscriptstyle HKR}$ n'envoie
pas le sous-espace $\gsI$ de $\gs$ dans le sous-espace $\GsI$ de $\Gs$:
en utilisant la d\'{e}composition
\begin{equation}
 \gsI ~~=~~ \Al\otimes \Lambda E_2\otimes \Lambda^+ E_1~~~\oplus~~~
         \Il \otimes \Lambda E_2
\end{equation}
on voit que l'image de $\psi_{\scriptscriptstyle HKR}$ contiendrait des
\'{e}l\'{e}ments
provenant du premier terme de la somme ci-dessus dont le facteur le plus
\`{a} droite est dans $E_2$ sans que le coefficient soit dans l'id\'{e}al, et un
tel op\'{e}rateur multidiff\'{e}rentiel ne serait plus dans $\GsI$.

 Il faut alors modifier les applications $\alpha^k$: il est bien connu que
 la multiplication ext\'{e}rieure
induit un isomorphisme d'espaces vectoriels
\begin{equation}
    \Phi:\Lambda E_2\otimes \Lambda E_1 \ra \Lambda (E_1\oplus E_2) =\Lambda E:
      T_2\otimes T_1 \mapsto T_2\wedge T_1.
\end{equation}
Soit $\iota^{(l,k-l)}:E_2^{\otimes l}\otimes E_1^{\otimes (k-l)}\ra
E^{\otimes k}$ l'injection induite par des sous-espaces $E_1,E_2$.
L'application $\hat{\alpha}^k:=\sum_{l=0}^k\iota^{(l,k-l)}\alpha^l\otimes
\alpha^{k-l}$ envoie $\big(\Lambda E_2\otimes \Lambda E_1\big)_k$ dans
$E^{\otimes k}$. Nous d\'{e}finissons l'application HKR modifi\'{e}e par
\begin{equation}
  {\psi^1}^k_{\scriptscriptstyle HKR}:\gs^k\ra \Gs^k: f\otimes T \mapsto
     f\otimes
     \hat{\alpha}^k\big(\Phi^{-1}(T)\big),
\end{equation}
et $\psi^1_{\scriptscriptstyle HKR}$ comme \'{e}tant la somme des tous les
${\psi^1}^k_{\scriptscriptstyle HKR}$. Ainsi le facteur le plus \`{a} droite
est toujours dans $E_1$ dans le cas o\`{u} son coefficient est dans $\Al$, et
n'est dans $E_2$ que si son coefficient est dans $\Il$. Par cons\'{e}quent,
$\psi^1_{\scriptscriptstyle HKR}$ envoie bien $\gsI$ dans $\GsI$.
En \'{e}crivant $\psi^{[1]}$ pour la restriction de
$\psi^1_{\scriptscriptstyle HKR}$
\`{a} $\gsI$ on a l'analogue suivant du
th\'{e}or\`{e}me \ref{applHKRusuel}:
\begin{theorem}\label{applHKR}
 Quel que soit l'entier positif $k$, les applications
 ${\psi^{[1]}}^k$ et
 $\pi^k_{\scriptscriptstyle HKR}$ ont les propri\'{e}t\'{e}s suivantes:
 \begin{enumerate}
  \item $\pi^k_{\scriptscriptstyle HKR}{\psi^{[1]}}^k =
          \mathrm{id}_{\gsI^k}$.
  \item $\mathsf{b}{\psi^{[1]}}^k=0$.
  \item Soit $\phi\in\GsI^k$. Si $\mathsf{b}\phi=0$ et
      $\pi^k_{\scriptscriptstyle HKR}\phi=0$, alors $\phi$ est un cobord,
      i.e. il existe $\psi\in\GsI^{k-1}$ tel que $\phi=\mathsf{b}\psi$.
  \item Soient $X,Y\in\gsI$, alors il existe $\phi\in \GsI$ tel que
      $[\psi^{[1]}X,\psi^{[1]}Y]_G
      -\psi^{[1]}[X,Y]_S=\mathsf{b}\phi$.
 \end{enumerate}
\end{theorem}
\begin{prooof}
 1. Le premier \'{e}nonc\'{e} est \'{e}vident. \\
 2. Puisque l'image de $\psi^{[1]}$ consiste en
 des op\'{e}rateurs multidiff\'{e}rentiels $1$-diff\'{e}rentiels, cet espace
 ne contient que des cocycles, d'o\`{u} le deuxi\`{e}me \'{e}nonc\'{e}.\\
 3. D'apr\`{e}s le troisi\`{e}me \'{e}nonc\'{e} du th\'{e}or\`{e}me \ref{applHKRusuel}, $\phi$
 est un cobord dans $\Gs$, mais puisque le morphisme connectant
 de la suite exacte longue correpondant \`{a} la suite exacte courte
 $\GsI\ra \Gs \ra \GsP$ s'annule, il s'ensuit que $\phi$ est aussi un
 cobord dans $\GsI$.\\
 4. D'apr\`{e}s 1. et le premier \'{e}nonc\'{e} de \ref{applHKRusuel},
 on trouve $\phi_1,\phi_2,\phi'\in\Gs$ tels que
 $\psi^{[1]}X=\psi_{\scriptscriptstyle HKR}X+\mathsf{b}\phi_1$,
 $\psi^{[1]}Y=\psi_{\scriptscriptstyle HKR}Y+\mathsf{b}\phi_2$ et
 $\psi^{[1]}[X,Y]=\psi_{\scriptscriptstyle HKR}[X,Y]+\mathsf{b}\phi'$.
 A l'aide de
 l'\'{e}nonc\'{e} 4. de \ref{applHKRusuel}, on voit que le membre de gauche de 4.
 est un \'{e}l\'{e}ment de $\GsI$ et un cobord dans $\Gs$, donc --\`{a} l'aide du
 m\^{e}me argument que pour 3.-- un cobord dans $\GsI$.
\end{prooof}

\section{Structure $G_\infty$ sur $\GsI$}

Dans cette section, nous allons tout d'abord rappeler les
d\'{e}finitions
des structures et morphismes $L_\infty$ et $G_\infty$. Puis, dans
la deuxi\`{e}me sous-section, nous rappellerons le
plan de la d\'{e}monstration de la conjecture de Deligne par
Tamarkin (il
existe une stucture $G_\infty$ sur $\Gs$) que nous adapterons au
cas o\`{u} l'espace est $\GsI$ dans la troisi\`{e}me sous-section pour
montrer l'existence d'une structure $G_\infty$ sur $\GsI$.

\subsection{Rappels et notations}
\noindent

\noindent Soit $\g$ un espace gradu\'{e}, nous noterons
$\L^{\cdot}\g$ la cog\`{e}bre colibre cocommutative
sur $\g$ et
$\L^{\cdot}{\underline{\g^{\otimes}}}$ l'espace gradu\'{e}
$$\somme_{m\geq 1, r_1+\cdots+r_n=m}{\underline{\g^{\otimes r_1}}}\w\cdots\w
{\underline{\g^{\otimes r_n}}}$$
o\`{u}  ${\underline{\g^{\otimes l}}}$ d\'{e}signe le quotient
de $\g^{\otimes l}$ par l'image des shuffles
d'ordre $l$
(pour plus de d\'{e}tails, voir \cite{Lo} ou \cite{Sta} par exemple).
On utlise la graduation suivante sur
$\L^{\cdot}{\underline{\g^{\otimes}}}$:
pour des \'{e}l\'{e}ments homog\`{e}nes $x_1^1,\dots,x_n^{r_n}$
dans $\g$ des degr\'{e}s respectivement $|x_1^1|,\dots,|x_n^{r_n}|$,
le degr\'{e} de
$x=(x_1^1\otimes \cdots\otimes
x_1^{r_1})\w \cdots\w (x_n^1\otimes \cdots\otimes x_n^{r_n})$
est
$$|x|=\sum_{i=1}^{r_1}|x_1^{i}|+\cdots+\sum_{i=1}^{r_n}|x_n^{i}|-n.$$
Il est bien connu que $\L^{\cdot}{\underline{\g^{\otimes}}}$
est une cog\`{e}bre colibre.
\begin{defi}
{~}

\begin{itemize}
\item Une structure d'alg\`{e}bre de Gerstenhaber \`{a} homotopie
pr\`{e}s (ou alg\`{e}bre $G_{\infty}$) sur un espace vectoriel
$\g$ est donn\'{e}e par une famille d'applications de degr\'{e} $1$~:
$$m^{r_1,\dots,r_n}:\g^{\otimes r_1}\w\cdots\w \g^{\otimes r_n}\rightarrow \g$$
telles que leur extension canonique \`{a} $\L^{.}{\underline{\g^{\otimes}}}$ satisfait
$$d\circ d=0$$
o\`{u} $$d=\sum_{m\geq 1} \sum_{r_1+\cdots+r_n=m} m^{r_1,\dots,r_n}.$$
\item
Une structure d'alg\`{e}bre de Lie \`{a} homotopie pr\`{e}s
(ou alg\`{e}bre $L_{\infty}$) sur un espace vectoriel
$\g$ est une alg\`{e}bre $G_\infty$ o\`{u} les applications
$m^{r_1,\dots,r_n}\!\!:\g^{\otimes r_1}\!\w\cdots\w \g^{\otimes r_n}\rightarrow \g$
sont nulles pour $r_i >1$.
\end{itemize}
\end{defi}
\noindent Consid\'{e}rons l'espace
$\gs=\Gamma(X,\ve^\cdot TX)$,
la structure d'alg\`{e}bre de Lie gradu\'{e}e de $\gs$ donn\'{e}e par
le crouchet de Schouten
$[-,-]_S$ peut \^{e}tre retraduite au moyen d'une application
$$m^{1,1}:\L^2\gs\rightarrow \gs$$
et la structure d'alg\`{e}bre commutative gradu\'{e}e, donn\'{e}e par le
produit
ext\'{e}ri\-eur
$\wedge$, au moyen d'une application
$$m^2:\gs^{\otimes 2}\rightarrow \gs.$$
Ces applications peuvent \^{e}tre naturellement \'{e}tendues en des
applications, toujours not\'{e}es $m^{1,1}$ and $m^2$, sur
$\L^{.}{\underline{\gs^{\otimes}}}.$
On v\'{e}rifie ais\'{e}ment que les applications
$m^{1,1}$ et $d=m^{1,1}+m^2$ v\'{e}rifient $$m^{1,1} \circ
m^{1,1}=0 \hbox{ et } d\circ d=0$$ de telle sorte que
$(\gs,m^{1,1})$ (resp. $(\gs, d)$) peut \^{e}tre vue comme une
alg\`{e}bre $L_\infty$ (resp. $G_{\infty}$).
Plus g\'{e}n\'{e}ralement, toute alg\`{e}bre de Gerstenhaber $(G,\mu,
[\,,\,])$
a une structure d'alg\`{e}bre $G_{\infty}$ canonique donn\'{e}e par
$m^{1,1}=[\,,\,]$, $m^2=\mu$, les autres applications \'{e}tant pos\'{e}es
nulles.

\smallskip

\noindent Consid\'{e}rons maintenant le complexe de Hochschild
${\Gs}=C^{.}(\Al,\Al)$ o\`{u} $\Al=C^{\infty}(X)$, vu comme
un espace vectoriel gradu\'{e}~:
${\Gs}=\somme_{k} C^{k}(\Al,\Al)$
(un \'{e}l\'{e}ment de $C^k(\Al,\Al)$ est de degr\'{e} $k-1$).
\'Equip\'{e} de crochet de Gerstenhaber $[-,-]_G$ et de la
diff\'{e}rentielle de Hochschild
$\mathsf{b}$, c'est une alg\`{e}bre de Lie diff\'{e}rentielle gradu\'{e}e (est
donc une alg\`{e}bre
$L_\infty$ avec $M^1=\mathsf{b}$ et $M^{1,1}=[-,-]_G$).
Deligne(\cite{Del}) a conjectur\'{e} que cet espace $\Gs$
peut aussi \^{e}tre muni d'une structure d'alg\`{e}bre
$G_{\infty}$ o\`{u} $M^2$ correspondrait au produit commutatif gradu\'{e}
usuel des
cocha\^{\i}nes~:
pour $E,F \in \Gs$ et
$x_1,\dots,x_{|E|+|F|+2}\in A$,
\begin{multline*}
M^2(E,F)(x_1,{\dots},x_{|E|+|F|+2})\cr
=(-1)^{\gamma}
E(x_1,{\dots},x_{|E|+1})F(x_{|E|+2},{\dots},x_{|E|+|F|+2})
\end{multline*}
 o\`{u}
$\gamma=(|F|+1)(|E|+1)$.
Dans la sous-section suivante, nous rappellerons le plan de la preuve
de cette conjecture donn\'{e}e par Tamarkin \cite{Tam}, et qui
utilise le foncteur de quantification-d\'{e}quantification
d'Etingof-Kazhdan \cite{EK}. Nous montrerons dans la troisi\`{e}me
sous-section comment adapter cette preuve pour montrer l'existence
d'une structure $G_\infty$ sur $\GsI$.
\begin{defi}
Un morphisme $L_\infty$ (respectivement $G_{\infty}$) entre deux alg\`{e}\-bres
$L_\infty$ (respectivement $G_{\infty}$) $(\g_1,d_1)$ et
$(\g_2,d_2)$ est un morphisme
$$\psi: (\L^{.}{\underline{\g_1^{\otimes}}},d_1)\rightarrow
(\L^{.}{\underline{\g_2^{\otimes}}},d_2)$$
 de cog\`{e}bre codiff\'{e}rentielle.
\end{defi}
\noindent Un morphisme $G_{\infty}$ $\psi$ entre deux alg\`{e}bres
$G_{\infty}$
$(\g_1, d_1)$ et $(\g_2, d_2)$ est donn\'{e} par une famille d'applications
$\psi^{[n]}$ o\`{u}
$$\psi^{[n]}=\sum_{r_1+\cdots+r_k=n} \psi^{r_1,\dots,r_k}$$
avec $\psi^{r_1,\dots,r_k}:\g_1^{\otimes r_1}\w\cdots\w \g_1^{\otimes r_k}\rightarrow
\g_2$ et qui satisfont, pour tout
$n$ l'\'{e}quation~:
$$\psi^{[\leq n]}\circ d_1^{[\leq n]}=d_2^{[\leq n]}\circ \psi^{[\leq n]}.$$
Ici, nous avons not\'{e}
$$\psi^{[\leq n]}=\sum_{r_1+\cdots+r_k\leq n} \psi^{r_1,\dots,r_k},$$
$$d_i^{[\leq n]}=\sum_{r_1+\cdots+r_k\leq n} d_{i}^{r_1,\dots,r_k}\quad (i=1,2)$$
o\`{u} $d_i^{r_1,\dots,r_k}$ sont les composantes
des codiff\'{e}rentielles $d_i$ ($i=1,2$).

\medskip

\subsection{Plan de la preuve de la conjecture de Deligne}
\noindent
Nous allons rappeler le plan de la preuve de la conjecture
de Deligne par Tamarkin telle qu'elle est
d\'{e}crite dans \cite{GH}. Notre but
est de construire une structure $G_\infty$, c'est-\`{a}-dire, une
diff\'{e}rentielle $D$ sur $\Gs$, puis sur le sous-espace $\GsI$,
satisfaisant toute
les deux, si
$$M=M^{1}+M^{1,1}+M^2+\cdots +
M^{p_1,\dots,p_n}+\cdots,$$
\begin{enumerate}
\item $M^1$ est le cobord de Hochschild $\mathsf{b}$ et $M^{1,1}$ est le
crochet $[-,-]_G$.
\item $M \circ M=0$.
\end{enumerate}
\noindent Dans la suite de cette section, $\h$ designera
$\Gs$ ou
sa sous-alg\`{e}bre $\GsI$. Le probl\`{e}me peut se reformuler comme suit~:
Soit $L=\oplus ~\underline{\h^{\otimes n}}$
la cog\`{e}bre de Lie colibre sur $\h$.
Puisque $L$ est une cog\`{e}bre de Lie colibre, une structure de
big\`{e}bre de Lie diff\'{e}rentielle sur
$L$ est donn\'{e}e par des applications de degr\'{e} $1$,
$l^n$~: $\underline{\h^{\otimes n}}\rightarrow \h$, correspondant
\`{a} la diff\'{e}rentielle, et par des applications
$l^{p_1,p_2}$~: $\underline{\h^{\otimes p_1}}\wedge
\underline{\h^{\otimes p_2}}
\rightarrow \h$, correspondant au crochet de Lie.
Ces applications s'\'{e}tendent de mani\`{e}re unique en des
d\'{e}rivations
de cog\`{e}bre
$L\to L$ et en un morphisme de cog\`{e}bre $L\wedge L\to L$ (toujours
not\'{e}s $l^m$ et $l^{p,q}$).
On a imm\'{e}diatement~:
\begin{lemma}\label{Lemma 1.1}
Supposons donn\'{e}e une structure de big\`{e}bre de Lie
diff\'{e}rentielle sur la cog\`{e}bre de
Lie $L$,
dont la diff\'{e}rentielle et le crochet de Lie sont determin\'{e}s
respectivement par les applications $l^n$
et
$l^{p_1,p_2}$ comme ci-dessus.
Alors $\h$ a une structure $G_\infty$ donn\'{e}e, pour tous $p,q,n\geq 1$, par
$$M^n=l^n, \qquad M^{p,q}=l^{p,q} \quad \mbox{ et }
\quad M^{p_1,\dots, p_r}=0 \mbox{ pour  } r\geq 3.$$
\end{lemma}
\begin{prooof}
Voir \cite{GH} pour la preuve d\'{e}taill\'{e}e.
\end{prooof}
Ainsi, pour obtenir la structure $G_\infty$ d\'{e}sir\'{e}e sur $\h$, il
est suffisant de d\'{e}finir une structure de big\`{e}bre de Lie
diff\'{e}rentielle sur $L$ donn\'{e}e par des applications $l^n$
et
$l^{p_1,p_2}$ avec $l^1=b$ et
$l^{1,1}=[-,-]_G$.

\medskip

Continuons \`{a} donner une formulation \'{e}quivalente de notre
probl\`{e}me~:
\begin{prop}\label{Theorem 1.2}
Supposons donn\'{e}e une structure de big\`{e}bre diff\'{e}rentielle
sur la cog\`{e}bre tensorielle colibre $T=\oplus_{n\geq 0}~\h^{\otimes
  n}$
dont la diff\'{e}rentielle et la
multiplication sont donn\'{e}es respectivement par des applications
$a^n$~: $\h^{\otimes n}\rightarrow \h$ et
$a^{p_1,p_2}$~: $\h^{\otimes p_1}\otimes \h^{\otimes p_2}
\rightarrow \h$. Alors on a une structure
de big\`{e}bre de Lie diff\'{e}rentielle sur la cog\`{e}bre de Lie
$L=\oplus_{n\geq 0}~\underline{\h^{\otimes n}}$,
dont la diff\'{e}rentielle et le crochet de Lie sont donn\'{e}s
respectivement
par les applications
$l^n$
et
$l^{p_1,p_2}$ o\`{u} $l^1=a^1$ et
$l^{1,1}$ est l'anti-symmetris\'{e}e de $a^{1,1}$.
\end{prop}
\begin{prooof}
La preuve utilise le th\'{e}or\`{e}me de quantification-d\'{e}quantifi\-cation
d'Etingof-Kazhdan (\cite{EK}). Elle est faite dans
\cite{Tam} et dans \cite{GH}.
\end{prooof}

\medskip

Ainsi, d\'{e}finir une structure de big\`{e}bre de Lie diff\'{e}rentielle sur $L$
donn\'{e}e par des applications $l^n$ et
$l^{p_1,p_2}$ avec $l^1=b$ et
$l^{1,1}=[-,-]_G$, est \'{e}quivalent \`{a} d\'{e}finir
une structure de big\`{e}bre diff\'{e}rentielle sur $T$
donn\'{e}e par des applications
$a^{n}$~: $\h^{\otimes n} \rightarrow \h$  et
$a^{p_1,p_2}$~: $\h^{\otimes p_1} \otimes \h^{\otimes p_2} \rightarrow
\h$
o\`{u} $a^1=b$
et $a^{1,1}$ est le produit
$\{-|-\}$ d\'{e}fini par
$$\{E|F\}(x_1,{\dots},x_{e+f-1})=\sum_{i\geq
0}{(-1)}^{|F|{\cdot}i}E(x_1,{\dots},x_i,F(x_{i+1},{\dots},
x_{i+f}),{\dots}),$$
pour $E,F$ dans $\h$
(il est clair, par d\'{e}finition, que l'anti-symm\'{e}tris\'{e} de
$\{-|-\}$
est le crochet de Gerstenhaber $[-,-]_G$).

\medskip

Cette deni\`{e}re  structure peut \^{e}tre obtenue, dans le cas o\`{u}
$\h=\Gs$
en utilisant le op\'{e}ration ``brace'' (d\'{e}finies dans \cite{Kad02} et
\cite{GV})
agissant sur le complexe de cocha\^{\i}nes de  Hochschild
$\Gs=C(\Al,\Al)$.
Dans la sous-section suivante, nous montrerons que
ces op\'{e}rations peuvent encore \^{e}tre utilis\'{e}es
dans le cas $\h=\GsI$, prouvant ainsi l'existence
d'une structure d'alg\`{e}bre $g_\infty$ sur $\GsI$.

\subsection{Existence d'une structure $G_\infty$ sur $\GsI$}
\noindent Comme nous l'avons vu pr\'{e}c\'{e}demment, pour
prouver l'existence de la structur d'alg\`{e}bre $G_\infty$
d\'{e}sir\'{e}e sur $\GsI$, nous devons
construire une
structure de  big\`{e}bre diff\'{e}rentielle sur $T
=\oplus_{n\geq 0}~\GsI^{\otimes
  n}$
donn\'{e}e par des applications
$a^{n}$~: $\GsI^{\otimes n} \rightarrow \GsI$  et
$a^{p_1,p_2}$~: $\GsI^{\otimes p_1} \otimes \GsI^{\otimes p_2} \rightarrow
\GsI$
o\`{u} $a^1=\mathsf{b}$
et $a^{1,1}$ est le produit
$\{-|-\}$ d\'{e}fini dans la section pr\'{e}c\'{e}dente. Pour cela nous
allons encore utiliser les op\'{e}rations braces qui se restreignent
\`{a} $\GsI$  puisque $\GsI$ v\'{e}rifie les propri\'{e}t\'{e}s de la
proposition \ref{PGsI}. Ces op\'{e}rations sont des applications~:
$a^{1,p}$~:
$\GsI \otimes \GsI^{\otimes p}\rightarrow \GsI$ ($p\geq 1$)
d\'{e}finies pour toutes cocha\^{\i}nes homog\`{e}nes $
E,F_1,\dots,F_p\in  \GsI^{\otimes p+1}$ et
$x_1,\dots, x_e \in A$ (avec $e= |E|+|F_1|+\dots +|F_p|+1$), par
\begin{multline*}
a^{1,p}(E_1\otimes (F_1 \otimes \dots \otimes F_p))(x_1 \otimes \cdots \otimes x_e)\cr
=\sum (-1)^{\tau} E(x_1,{\dots},x_{i_1},F_1(x_{i_1+1},{\dots}),
{\dots},F_p(x_{i_p+1},{\dots}),{\dots})
\end{multline*}
o\`{u} $\tau=\sum_{k=1}^pi_k(|F_k|+1)$.
Les applications $a^{1,p}:\GsI \otimes \GsI^{\otimes p}\rightarrow \GsI$
et $a^{q\geq 2,p}=0$ donnent une unique structure
de big\`{e}bre sur l'alg\`{e}bre cotensorielle colibre
$T=\oplus_{n\geq 0}\GsI^{\otimes n}$.
De m\^{e}me, en posant $a^1$ comme \'{e}tant \'{e}gale au cobord de
Hochschild $\mathsf{b}$, $a^2$ le
product $M^2$, et $a^{q\geq 3}=0$, on obtient une unique
structure de big\`{e}bre diff\'{e}rentielle
sur la cog\`{e}bre tensorielle $T$.
Le Th\'{e}or\`{e}me 3.1 dans \cite{GV} nous assure
que ces applications induisent une structure de big\`{e}bre
diff\'{e}rentielle sur la cog\`{e}bre colibre $T$, ce qui nous donne
la structure d'alg\`{e}bre $G_\infty$ voulue sur $\GsI$.
Par construction, les applications $M^{p_1,\dots,p_k}$ sont nulles
pour $k >2$.
De plus, l'application $M^2$ co\"{\i}ncide,
\`{a} un cobord de Hochschild pr\`{e}s,
avec le produit $M^2$ car,
apr\`{e}s passage \`{a} la cohomologie, elles donnent toutes les
deux la m\^{e}me application $m^2$, correspondant au produit
ext\'{e}rieur $\wedge$ de l'alg\`{e}bre
de Gerstenhaber $(\gs,[-,-]_S,\wedge)$.

\section{Obstructions \`{a} la formalit\'{e} entre $\gsI$ et $\GsI$}

Dans cette section, nous allons reprendre le plan de la d\'{e}monstration
de Tamarkin (\cite{Tam}, et aussi \cite{GH}) du th\'{e}or\`{e}me de formalit\'{e}
$G_\infty$ et montrer que les
obstructions \`{a} la construction d'un morphisme $G_\infty$ entre $(\gsI,d)$
et $(\GsI,D)$ (o\`{u} $d$ et $D$ d\'{e}finissent les structures d'alg\`{e}bre
$G_\infty$ sur $\gsI$ et $\GsI$)
se trouvent dans le groupe de cohomologie de
$\left(\!{\rm Hom}({\ve^\cdot\underline{\gsI^{\otimes \cdot}}},\gsI),
\left[[-,-]_S+\w,-\right]\right)$. Dans la premi\`{e}re sous-section,
nous montrerons qu'il existe une autre structure d'alg\`{e}bre $G_\infty$
donn\'{e}e par une diff\'{e}rentielle $d'$ et un morphisme
$G_\infty$, $\psi$, entre $(\gsI,d')$ et $(\GsI,D)$.
Dans la deuxi\`{e}me sous-section nous montrerons que
si le groupe de cohomologie du complexe
$\left(\!{\rm Hom}({\ve^\cdot\underline{\gsI^{\otimes \cdot}}},\gsI),
\left[[-,-]_S+\w,-\right]\right)$ est trivial, on peut
construire un morphisme $G_\infty$, $\psi'$, entre
$(\gsI,d)$ et $(\gsI,d')$. La compos\'{e}e $\psi \circ \psi'$ sera donc alors
bien un morphisme $G_\infty$ entre $(\gsI,d)$ et $(\GsI,D)$.

\subsection{Morphisme $G_\infty$ entre $(\gsI,d')$ et $(\GsI,D)$}

\noindent
Dans cette sous-section, nous allons prouver la proposition suivante~
\begin{prop}\label{Theorem 2.1}
Il existe une structure $G_\infty$ donn\'{e}e par une
diff\'{e}rentielle $d'$ sur
$\gu $ et un morphisme $G_\infty$
$\psi$~: $(\gu,d')\rightarrow (\gd,D)$
tels que la restriction $\psi^1:\gsI\to \GsI$  est
l'application de Hochschild-Kostant-Rosenberg d\'{e}finie en section 2
(voir theor\`{e}me \ref{theoHKR} et th\'{e}or\`{e}me \ref{applHKR}).
\end{prop}
\begin{prooof}
Pour $n\geq 0$, on notera $$\gsI^{[n]}=\bigoplus_{p_1+\cdots +p_k=n}
~\underline{\gsI^{\otimes p_1}}\wedge
\cdots \wedge  \underline{\gsI^{\otimes p_k}}$$ et
$\gsI^{[\leq n]}=\sum_{k \leq n}\gsI^{[k]}$. De m\^{e}me, on notera
$\GsI^{[n]}=\bigoplus_{p_1+\cdots +p_k=n}
~\underline{\GsI^{\otimes p_1}}\wedge
\cdots \wedge  \underline{\GsI^{\otimes p_k}}$ et
$\GsI^{[\leq n]}=\sum_{k \leq n}\GsI^{[k]}$.
Soit $D^{p_1,\dots,p_k}$~: $\underline{\GsI^{\otimes p_1}}\wedge
\cdots \wedge \underline{\GsI^{\otimes p_k}}
\rightarrow \GsI$ les composantes de la diff\'{e}rentielle $D$ definissant la
structure
$G_\infty$ de $\GsI$.
On notera  $D^{[n]}$ et $D^{[\leq n]}$ les sommes
$$D^{[n]}=\sum_{p_1+\cdots+p_k=n}D^{p_1,\dots,p_k}
\qquad \mbox{et} \qquad D^{[\leq n]}=\sum_{p\leq n} D^{[p]}.$$
Clairement,  $D=\sum_{n\geq 1} D^{[n]}$.
De la m\^{e}me mani\`{e}re, on note
$${d'}_1^{[ n]}=\sum_{p_1+\cdots+p_k=n}d^{'p_1,\dots,p_k} \qquad \mbox{et} \qquad
{d'}_1^{[\leq n]}=\sum_{1 \leq k \leq n} {d'}_1^{[k]}.$$
D'apr\`{e}s la section pr\'{e}c\'{e}dente, un morphisme $\psi$~:
$(\gu,d')\rightarrow
(\gd,D)$ est d\'{e}ter\-min\'{e} de mani\`{e}re unique par
ses composantes $\psi^{p_1,\dots,p_k}:
\underline{\gsI^{\otimes p_1}}\wedge \cdots
\wedge \underline{\gsI^{\otimes p_k}}\rightarrow \GsI$.
On pose alors encore
$$\psi^{[n]}=\sum_{p_1+\cdots+p_k=n}\psi^{p_1,\dots,p_k}\qquad
\mbox{et} \qquad
\psi^{[\leq n]}=\sum_{1 \leq k \leq n} \psi^{[k]},$$
et aussi $d'=\sum_{n \geq 1} {d'}^{[n]}$ et $\psi=\sum_{n \geq 1} \psi^{[n]}$.

\medskip

Pour construire la diff\'{e}rentielle $d'$ et le morphisme
$\psi$,
on va construire les applications ${d'}^{[n]}$ et
$\psi^{[n]}$ par r\'{e}currence. Pour les premiers termes,
on pose
${d'}^{[1]}=0$ et $\psi^{[1]}$ est l'isomorphisme de
Hochschild-Kostant-Rosenberg du theor\`{e}me \ref{applHKR}.

\medskip

Supposons construites les applications
$({d'}^{[i]})_{i \leq n-1}$ et $(\psi^{[i]})_{i
\leq n-1}$ v\'{e}rifiant les conditions
$$\psi^{[\leq n-1]} \circ {d'}^{[\leq n-1]}=D^{[\leq n-1]} \circ \
\psi^{[\leq n-1]}$$ sur $\gsI^{[\leq n-1]}$
et
$${d'}^{[\leq n-1]} \circ {d'}^{[\leq n-1]}=0$$
sur $\gsI^{[\leq n]}$.
Ces conditions suffisent pour nous assurer que $d'_1$
est une diff\'{e}ren\-tielle et que $\psi$ est un morphisme de
cog\`{e}bre diff\'{e}rentielle.
Si l'on reformule l'identit\'{e} $\psi \circ {d'}=D \circ \psi$
sur $\gsI^{[n]}$, on obtient
\begin{equation}
\psi^{[\leq n]} \circ {d'}^{[\leq n]}=D^{[\leq n]}\circ\psi^{[\leq
n]}.
\label{E}
\end{equation}
\noindent Si l'on tient compte, maintenant, du fait que
${d'}^{[1]}=0$, et que sur $\gsI^{[n]}$ on a
$\psi^{[k]} \circ {d'}^{[l]}=D^{[k]}\circ\psi^{[l]}=0$ pour
$k+l>n+1$, l'identit\'{e} (\ref{E}) devient
\begin{equation}
\psi^{[1]}{d'}^{[n]} + B = {D}^{[1]} \psi^{[n]} + A
\label{E'}
\end{equation}
\noindent o\`{u} $B=\sum_{k=2}^{n-1}\psi^{[\leq n-k+1]}{d'}^{[k]}$
et $A={D}^{[1]}\psi^{[\leq n-1]}+\sum_{k=2}^{n}{D}^{[k]}\psi^{[\leq n-k+1]}$
(on oublira, par la suite, le symb\^{o}le de composition $\circ$).
Le terme ${D}^{[1]}$ dans (\ref{E'}) est le cobord de Hochschild
$b$. Alors, gr\^ace au th\'{e}or\`{e}me de Hochschild-Kostant-Rosenberg
entre $\gsI$ et $\GsI$,
l'identit\'{e} (\ref{E'}) est \'{e}quivalente au fait que la
cocha\^{\i}ne $B-A$  est un
cocycle de Hochschild. Ainsi, pour montrer l'existence de ${d'}^{[n]}$
et $\psi^{[n]}$, il sera suffisant de montrer que
\begin{equation}
{D}^{[1]}(B-A)=0
\label{E1}
\end{equation}
\noindent et de montrer ensuite que pour n'importe quel choix
de cobord $\psi^{[n]}$, on aura toujours
\begin{equation}
{d'}^{[\leq n]}{d'}^{[\leq n]}=0\hbox{ sur }
\gsI^{[\leq n+1]}.
\label{E2}
\end{equation}
$\bullet$ Commen\c cons par construire ${d'}^{[2]}$:
pour $n=2$, on obtient $A={D}^{[1]}\psi^{[1]}+
{D}^{[2]}\psi^{[1]}$ et $B=0$ et donc
$$\psi^{[1]}{d'}^{[2]} = {D}^{[1]}( \psi^{[2]} + \psi^{[1]})
+{D}^{[2]}\psi^{[1]}.$$
Ainsi ${d'}^{[2]}$ est
l'image de ${D}^{[2]}$ par la projection sur
la cohomologie de $\GsI$ et puisque l'application
de Hochschild-Kostant-Rosenberg $\psi^{[1]}$
est injective de $\gsI$ ($=H(\GsI,b={D}^{[1]})$) vers $\GsI$, on obtient
$${d'}^{[2]}=d^{[2]}.$$
$\bullet$ Prouvons maintenant (\ref{E1}): on a
${D}^{[1]}(-A)=-\sum_{k=2}^{n}{D}^{[1]}{D}^{[k]}\psi^{[\leq n-k+1]}$.
En utilisant  $D\circ D=0$, on obtient
\begin{align*}{D}^{[1]}(-A)&=\sum_{k=2}^{n}\left(\sum_{l=2}^k D^{[l]}D^{[k+1-l]}\right)
\psi^{[\leq n-k+1]}\\
&=\sum_{l=2}^{n}D^{[l]}\left(
\sum_{k=l}^n D^{[k+1-l]}\psi^{[\leq n-k+1]}\right).\end{align*}
On a clairement
$\sum_{k=l}^n D^{[k+1-l]}\psi^{[\leq n-k+1]}=
\sum_{k=1}^{n-l+1} D^{[k]}\psi^{[\leq n-k+2-l]}$.
Si on utilise encore le fait que
$d^{[a]}d^{[b]}\psi^{[c]}=0$ sur $\gsI^{[n]}$ pour $a+b+c > n+2$, on
peut ajouter
des termes
$(\psi^{[n-k+2-l+k']})_{0 \leq k'\leq k-1}$ to $\psi^{[\leq n-k+2-l]}$
sans changer l'\'{e}galit\'{e} pr\'{e}c\'{e}dente. On a ainsi
$${D}^{[1]}(-A)=\sum_{l=2}^{n}D^{[l]}\left(
\sum_{k=1}^{n-l+1} D^{[k]}\right)\psi^{[\leq n+1-l]}
=\sum_{l=2}^{n}D^{[l]}
D^{[\leq n+1-l]}\psi^{[\leq n+1-l]}.$$
Puisque $\left(D^{[l]}\right)_{l \geq 2}$ envoie  $\GsI^{[\leq k]}$ sur
$\GsI^{[\leq k-1]}$, l'\'{e}galit\'{e} pr\'{e}c\'{e}dente n'a de termes non-triviaux
que sur $\gsI^{[\leq n-1]}$.
On peut alors appliquer l'hypoth\`{e}se de r\'{e}cur\-rence
$\psi^{[\leq k]} {d'}^{[\leq k]}=D^{[\leq k]}
\psi^{[\leq k]}$ sur $\gsI^{[\leq k]}$ pour $k \leq n-1$. On obtient
$${D}^{[1]}(-A)=\sum_{l=2}^{n}D^{[l]}
\psi^{[\leq n+1-l]}{d'}^{[\leq n+1-l]}.$$
On a maintenant
$${D}^{[1]}(B-A)={D}^{[1]}\sum_{k=2}^{n-1}\psi^{[\leq n-k+1]}{d'}^{[k]}
+\sum_{l=2}^{n}D^{[l]}
\psi^{[\leq n+1-l]}{d'}^{[\leq n+1-l]}.$$
Le terme correspondant \`{a} $l=n$
s'annule puisque ${d'}^{[1]}=0$.
En projetant toujours sur $\gsI^{[n]}$
les applications de type $D^{[a]}\psi^{[b]}{d'}^{[c]}$, on peut ajouter
des applications $\psi^{[p+p']}$ ($p'\geq 0$) \`{a} $\psi^{[\leq p]}$.
On obtient alors, apr\`{e}s r\'{e}indexation,
$${D}^{[1]}(B-A)=
\sum_{l=2}^{n-1}D^{[\leq n+1-l]}\,\psi^{[\leq n-1]} {d'}^{[l]}.$$
On a ainsi montr\'{e} que
${D}^{[1]}(B-A)=
\sum_{l=2}^{n-1}D^{[\leq n+1-l]}\,\psi^{[\leq n+1-l]} {d'}^{[l]}
$.
Puisque ${d'}^{[l]}\left(\gsI^{[\leq k]}\right)$ $\subset \gsI^{[\leq k-1]}$,
on peut encore utiliser l'hypoth\`{e}se de r\'{e}currence et on obtient
$${D}^{[1]}(B-A)=
\sum_{l=2}^{n-1}
\psi^{[\leq n+1-l]}{d'}^{[\leq n+1-l]}\,{d'}^{[ l]}=0$$
car ${d'}^{[1]}=0$ et
${d'}^{[\leq n-1]}{d'}^{[\leq n-1]}=0$ sur $\gsI^{[\leq n]}$,
encore par hypoth\`{e}se de r\'{e}currence

\noindent$\bullet$ On prouve finalement (\ref{E2}),
c'est-\`{a}-dire ${d'}^{[\leq n]}{d'}^{[\leq n]}=0$ sur
$\gsI^{[\leq n+1]}$. Puisque $\psi^{[1]}$ est un quasi-isomorphisme
entre $(\gsI,0)$ et $(\GsI,b=D^{[1]})$, il nous suffit de montrer que
$$\psi^{[1]}\;{d'}^{[\leq n]}\,{d'}^{[\leq n]}
\hbox{ est un bord sur }
\gsI^{[\leq n+1]}.$$
Apr\`{e}s avoir projet\'{e} les applications,
on obtient l'identit\'{e}  suivante sur $\gsI^{[\leq n+1]}$:
$$\psi^{[1]}\,{d'}^{[\leq n]}\,{d'}^{[\leq n]}=
\psi^{[\leq n]}\,{d'}^{[\leq n]}\,{d'}^{[\leq n]}.$$
D'apr\`{e}s la d\'{e}finition de ${d'}^{[\leq n]}$ on peut \'{e}crire
$\psi^{[\leq n]}\,{d'}^{[\leq n]}=
{D}^{[\leq n]}\,\psi^{[\leq n]}$ car ${d'}^{[\leq n]}$
envoie $\gsI^{[\leq n+1]}$ sur $\gsI^{[\leq n]}$.
Ainsi, il sera suffisant de montrer que
${D}^{[\leq n]}\psi^{[\leq n]}{d'}^{[\leq n]}$  est
un bord quand on le restreint \`{a} $\gsI^{[\leq n+1]}$.
On a alors
$${D}^{[\leq n]}\psi^{[\leq n]}{d'}^{[\leq n]}
=b\psi^{[\leq n]}{d'}^{[\leq n]}+ \sum_{2 \leq k \leq n}
{D}^{[\leq k]}\psi^{[\leq n]}{d'}^{[\leq n]}.$$
Puisque $\sum_{2 \leq k \leq n}
{D}^{[\leq k]}$ envoie $\GsI^{[\leq k]}$ sur
$\GsI^{[\leq k-1]}$, l'expression
$$\sum_{2 \leq k \leq n}
{D}^{[\leq k]}\psi^{[\leq n]}{d'}^{[\leq n]}$$
a des composantes non nulle seulement sur $\gsI^{[\leq n+1]}$.
Sur ce dernier espace, on a
$$\sum_{2 \leq k \leq n}
{D}^{[\leq k]}\,\psi^{[\leq n]}\,{d'}^{[\leq n]}=
\sum_{2 \leq k \leq n}
{D}^{[\leq k]}\,D^{[\leq n]}\,\psi^{[\leq n]},$$ d'apr\`{e}s la
d\'{e}finition de
${d'}^{[\leq n]}$.
Ainsi, on obtient les identit\'{e}s suivantes sur $\gsI^{[\leq n+1]}$:
$${D}^{[\leq n]}\,\psi^{[\leq n]}\,{d'}^{[\leq n]}=
b\psi^{[\leq n]}\,{d'}^{[\leq n]}
-b\,D^{[\leq n]}\,\psi^{[\leq n]}
+D^{[\leq n]}\,D^{[\leq n]}\,\psi^{[\leq n]}$$
$$=b\,\psi^{[\leq n]}\,{d'}^{[\leq n]}
-b\,D^{[\leq n]}\,\psi^{[\leq n]}$$
puisque $D\circ D=0$.
\end{prooof}

\subsection{Morphisme $G_\infty$ entre $(\gsI,d)$ et $(\gsI,d')$}

\noindent
Dans cette sous-section  nous allons montrer la proposition~:
\begin{prop}\label{Theorem 3.1}
Si le complexe $\left(\!\Hom(\gu,\gsI),[m_1^{1,1}\!\!\!+\!m_1^2,-]\!\right)$
est concen\-tr\'{e} en bidegr\'{e} $(0,0)$,
il existe un morphisme
$G_\infty$, $\psi '$~:
$(\gu,d) \rightarrow (\gu, d')$ tel que
la restriction ${\psi'}^{[1]}:\gsI\to \gsI$ est l'identit\'{e}.
\end{prop}

\smallskip

Nous utiliserons les m\^{e}mes notations pour
$\gsI^{[n]}$, $\gsI^{[\leq n]}$,
${d'}^{[n]}$ et ${d'}^{[\leq n]}$ que dans la sous-section pr\'{e}c\'{e}dente.
Nous noterons aussi
$$d=\sum_{n\geq 1} d^{[n]} \qquad \mbox{et} \qquad
{d}^{[\leq n]}=\sum_{1 \leq k \leq n} {d}^{[k]}$$ et similarly
$$\psi'=\sum_{n \geq 1} \psi^{[n]}\qquad \mbox{et} \qquad
{\psi'}^{[\leq n]}=\sum_{1 \leq k \leq n} {\psi'}^{[n]}.$$
\smallskip
\begin{prooof}
Nous allons construire les applications ${\psi'}^{[n]}$
par r\'{e}currence comme
dans la sous-section pr\'{e}c\'{e}dente.
Pour ${\psi'}^{[1]}$ nous poserons~:
$${\psi'}^{[1]}=\Id \hbox{ (l'application identit\'{e})}.$$
Supposons construites les applications $({\psi'}^{[i]})_{i\leq n-1}$
satisfaisant
$${\psi'}^{[\leq n-1]} {d}^{[\leq n]}={d'}^{[\leq n]}
{\psi'}^{[\leq n-1]}$$ sur $\gsI^{[\leq n]}$ (${d'}^{[\leq n]}$
envoie $\gsI^{[\leq l]}$ sur
$\gsI^{[\leq l-1]}$).
L'\'{e}quation ${\psi'}  {d}={d'}  \psi'$
sur $\gsI^{[n+1]}$, nous donne
\begin{equation}
{\psi'}^{[\leq n]}\, {d}^{[\leq n+1]}={d'}^{[\leq n+1]}\,{\psi'}^{[\leq
n]}.
\label{EE'}
\end{equation}
Puisque $d^{[i]}=0$ pour $i\not=2$,
${d'}^{[1]}=0$
et que sur $\gsI^{[n+1]}$ on a ${\psi'}^{[ k]} {d}^{[l]}=
{d'}^{[\leq k]}{\psi'}^{[l]}=0$ pour $k+l > n+2$, l'identit\'{e} (\ref{EE'})
devient
$${\psi'}^{[\leq n]}\,d^{[2]}=\sum_{k=2}^{n+1} {d'}^{[k]}\,{\psi'}^{[\leq
n-k+2]}.$$
Nous avons vu dans la sous-section pr\'{e}c\'{e}dente que
${d'}^{[2]}={d}^{[2]}$. Ainsi (\ref{EE'}) devient \'{e}quivalent \`{a}
$${d}^{[2]}{\psi'}^{[\leq n]}-{\psi'}^{[\leq n]}{d}^{[2]}=
\left[{d}^{[2]},{\psi'}^{[\leq n]}\right]=
-\sum_{k=3}^{n+1} {d'}^{[k]}{\psi'}^{[\leq
n-k+2]}.$$
Remarquons que  $d^{[2]}=m_1^{1,1}+m_1^2$.
Le complexe
$({\rm Hom}(\gu,\gsI),[d^{[2]},-])$ est acyclique,
ce qui nous assure que la construction de
${\psi'}^{[\leq n]}$ sera possible quand
$\sum_{k=3}^{n+1} {d'}^{[k]}{\psi'}^{[\leq
n-k+2]}$ est un cocycle dans ce complexe.
Ainsi, pour finir la preuve, il nous reste \`{a} montrer que
\begin{equation}
\left[d^{[2]},\sum_{k=3}^{n+1} {d'}^{[k]}\,{\psi'}^{[\leq n-k+2]}\right]=0
\hbox{ sur }\gsI^{[n+1]}.
\label{E1'}
\end{equation}
On a $$D_n=
\left[d^{[2]},\sum_{k=3}^{n+1} {d'}^{[k]}\,{\psi'}^{[\leq n-k+2]}\right]=
\left[d^{[2]},\sum_{k=1}^{n-1}{d'}^{[n+2-k]}\,{\psi'}^{[\leq k]}\right].$$
Il s'ensuit que l'on peut \'{e}crire
\begin{equation}-D_n=
\sum_{k=1}^{n-1}\left[d^{[2]},{d'}^{[n+2-k]}\right]{\psi'}^{[\leq k]}-
\sum_{k=1}^{n-1}{d'}^{[n+2-k]}\left[d^{[2]},{\psi'}^{[\leq k]}\right].\label{EEE1'}\end{equation}
En utilisant l'hypoth\`{e}se de r\'{e}currence pour
$\left({\psi'}^{[\leq k]}\right)_{k \leq n-1}$, on obtient
$$[{d}^{[2]},{\psi'}^{[\leq k]}]=
-\sum_{l=3}^{k+1} {d'}^{[l]}\,{\psi'}^{[\leq
k-l+2]}=
-\sum_{l=1}^{k-1}{d'}^{[k+2-l]}\,{\psi'}^{[\leq l]}$$ sur $\gsI^{[\leq k+1]}$.
L'\'{e}quation \ref{EEE1'} devient alors
$$-D_n=\sum_{k=1}^{n-1}\left[d^{[2]},{d'}^{[n+2-k]}\right]{\psi'}^{[\leq k]}
+\sum_{k=1}^{n-1}{d'}^{[n+2-k]}\left(
\sum_{l=1}^{k-1}{d'}^{[k+2-l]}\,{\psi'}^{[\leq l]}\right).$$
Finalement on a
$$-D_n=\sum_{k=1}^{n-1}\left[d^{[2]},{d'}^{[n+2-k]}\right]\,{\psi'}^{[\leq k]}
+\sum_{l=1}^{n-2}\left(\sum_{k=l+1}^{n-1}{d'}^{[n+2-k]}{d'}^{[k+2-l]}\right)
{\psi'}^{[\leq l]}.$$
Ceci implique que $$~-D_n=\sum_{k=1}^{n-1}\left(
\left[{d'}^{[2]},{d'}^{[n+2-k]}\right]+\sum_{p=k+1}^{n-1}
{d'}^{[n+2-p]}{d'}^{[p+2-k]}\right){\psi'}^{[\leq k]}.$$ Mais
les applications $$\left[{d'}^{[2]},{d'}^{[n+2-k]}\right]+\sum_{p=k+1}^{n-1}
{d'}^{[n+2-p]}\,{d'}^{[p+2-k]}=
\sum_{q=2}^{n+2-k}{d'}^{[q]}\,{d'}^{[n+4-q-k]}$$ sont nulles
car ${d'}\circ{d'}=0$ sur
$\gsI^{[\leq n+2-k]}$. Ceci donne le r\'{e}sultat.
\end{prooof}

\section{De la formalit\'{e} \`{a} la star-repr\'{e}sentation}

Dans cette section, nous consid\'{e}rons un champ
de tenseurs de Poisson $P$ compatible avec la
sous-vari\'{e}t\'{e} $C$.
Nous allons montrer comment l'existence d'un morphisme $G_\infty$
entre $\gsI$ et $\GsI$, construit dans la
section pr\'{e}c\'{e}dente (sous r\'{e}serve que le complexe $
\left(\!{\rm Hom}({\ve^\cdot\underline{\gsI^{\otimes \cdot}}},\gsI),
\left[[-,-]_S+\w,-\right]\right)$
est concentr\'{e} en bidegr\'{e} $(0,0)$), nous permet construire un star-produit
d\'{e}fini par des cocha\^{\i}nes compatibles et donc tel  que $\IlL$ est
un id\'{e}al \`{a} gauche dans
$(\AlL,*)$ (ce qui impliquera la star-repr\'{e}sentation).

\smallskip

Nous allons montrer dans un premier temps que le morphisme $G_\infty$
de la section pr\'{e}dente entre les alg\`{e}bres $G_\infty$, $(\gsI,d)$ et
$(\GsI,D)$,
se restreint en un morphisme $L_\infty$ entre $\gsI$ et $\GsI$ (vues cette fois
comme
alg\`{e}bre $L_\infty$). Remarquons, tout d'abord que les diff\'{e}rentielles $d$
et $D$ d\'{e}finissant les structures d'alg\`{e}bre $G_\infty$ sur $\gsI$ et $\GsI$
se restreignent sur $\ve^\cdot \gsI$ et $\ve^\cdot \GsI$ respectivement
en les codiff\'{e}rentielles $m^{1,1}$ ($=[-,-]_S$) et
$M^1+M^{1,1}$ ($=\mathsf{b}+[-,-]_G$) respectivement
(en effet, d'apr\'{e}s la fin de la Section 3,
les applications $M^{p_1,\dots,p_k}$ sont nulles
pour $k >2$). Ainsi, on restreignant l'application
$\phi=\psi  \circ \psi ': (\gu,d) \rightarrow (\gd,D)$
\`{a} $\phi$~: $\ve^\cdot\gsI \to  \ve^\cdot\GsI$, on obtient bien
un morphisme $L_\infty$ entre $(\gsI,m^{1,1})$ et $(\GsI,M^{1}+M^{1,1})$.

\bigskip

\`A partir de maintenant, on va supposer que $X$
est une vari\'{e}t\'{e} de Poisson
munie d'un champ de tenseur $P$ (satisfaisant $[P,P]_S=0$)
compatible avec  la sous-vari\'{e}t\'{e} $C$.
Soit $\hbar$ un param\`{e}tre formel et \'{e}tendons
$\phi$ en un morphisme $L_\infty$ $\real[[\hbar]]$-lin\'{e}aire~:
$\ve^\cdot\gsI[[\hbar]] \to  \ve^\cdot\GsI[[\hbar]]$.
Ce morphisme $L_\infty$ nous permet de
construire un star-produit $\star$ sur $X$ (see \cite{BFFLS78}),
``compatible'' (c'est-\`{a}-dire tel que $\IlL$ est
un id\'{e}al \`{a} gauche dans
$(\AlL,*)$)~: posons
$P_\hbar=\sum_{n\geq 0} \hbar^n\ve^n  P \in \ve^\cdot \gsI$
o\`{u} $\ve^n P=\underbrace{P \wedge \cdots \wedge P}_{n~\tim}$
(ici $\wedge$ n'est pas le produit ext\'{e}rieur des champs de vecteurs
mais $a \wedge b$ est
un \'{e}l\'{e}ment de $\gsI \wedge \gsI$).
D\'{e}finissons maintenant $m_\star=\phi(P_\hbar)$, on obtient
\begin{equation}[m_\star,m_\star]_G=0.\label{mstar}\end{equation}
Ceci est la cons\'{e}quence de la d\'{e}finition d'un morphisme
$L_\infty$ et du fait $[P,P]_S=0$ implique
$m^{1,1}(P_\hbar)=0$.
L'application $m_\star$ est un \'{e}l\'{e}ment de
$\GsI[[\hbar]]$ de degr\'{e} $1$,  elle d\'{e}finit donc une
application dans $C^2(\Al,\Al)[[\hbar]]$,
o\`{u} $C^2(\Al,\Al)[[\hbar]]$
d\'{e}signe l'ensemble de applications $\real[[\hbar]]$-bilin\'{e}aires
dans $C^2(\Al,\Al)$. L'identit\'{e} \ref{mstar} implique que  $m_\star$ est
un produit associatif sur $\Al[[\hbar]]$.
Enfin, par d\'{e}finition de $\phi$, on a:
$$m_\star=m+\hbar \phi^1(P) + \sum_{n\geq 2} \hbar^n \phi^n(P,\dots,P),$$
o\`{u} $\phi^1(P)$ est l'image du crochet de Poisson par le morphisme
de Hochschild-Kostant-Rosenberg d\'{e}fini Section 2. Notons que
$\phi^1(P)(f,g)-\phi^1(P)(g,f)=2\{f,g\}$ car $\phi^1$
est le morphisme de Hochschild-Kostant-Rosenberg usuel
auquel on a rajout\'{e} des termes sym\'{e}triques. Enfin
les $\phi^n(P,\dots,P)$ sont des cocha\^{\i}nes compatibles. Ceci nous donne le
r\'{e}sultat souhait\'{e}.

\section{Calcul des obstructions}

Dans les sections pr\'{e}c\'{e}dentes, nous avons vu que,
dans le cas o\`{u} $X=\real^n$, $C=\real^{n-\nu}$, les obstructions
\`{a} la construction de star-repr\'{e}sentations r\'{e}sidaient dans le
groupe de cohomologie de
$\left(\!{\rm Hom}({\ve^\cdot\underline{\gsI^{\otimes \cdot}}}, \gsI),
\left[[-,-]_S+\w,-\right]\right)$. Dans cette section, nous nous proposons
de calculer ces obstructions. plus pr\'{e}cis\'{e}\-ment,  on munit un \'{e}l\'{e}ment
$x\in
\underline{\gsI^{\otimes p_1}}\wedge \dots \underline{\gsI^{\otimes p_n}}$
du bidegr\'{e} $(\sum_{i=1}^np_i-1, n-1)$.
Cette graduation donne une structure de bicomplexe \`{a} l'espace vectoriel
$\left(\!{\rm Hom}({\ve^\cdot\underline{\gsI^{\otimes \cdot}}}, \gsI),
\left[[-,-]_S+\w,-\right]\right)$ pour laquelle $[-,-]_S$ est de bidegr\'{e}
$(0,1)$ et $\w$ de bidegr\'{e} $(1,0)$ {\it cf.}~\cite{Tam}.
Dans la premi\`{e}re sous-section,
nous montrerons que  le complexe
$$\left(\!{\rm Hom}({\ve^\cdot\underline{\gsI^{\otimes \cdot}}},\gsI),
\left[[-,-]_S+\w,-\right]\right)$$
est concentr\'{e} en bidegr\'{e} $(0,0)$
si le complexe
$$\left({\Hom}_{{\g}}(\ve_{{\g}}^\cdot
\g \ot_{{\gsI}_+}\Omega_{{\gsI}_+},\gsI),\delta^{1,1} \right),$$
est concentr\'{e} en degr\'{e} $0$,
o\`{u} $\delta=\delta^{1,1}+\delta^2$ est l'application
duale de $[m^{1,1}+m^2,-]=\left[[-,-]_S+\w,-\right]$
et $\Omega_{\gsI}$ est le module des $1$-formes diff\'{e}rentielles
de K\"ahler de l'alg\`{e}bre $\gsI$.
Enfin, dans la deuxi\`{e}me sous-section, nous discuterons les
cas o\`{u} nous conjecturons que
ce complexe est concentr\'{e} en degr\'{e} $0$.
On dira parfois abusivement qu'un complexe est acyclique pour dire qu'il est
concentr\'{e} en degr\'{e} (ou bidegr\'{e}) $0$ ($(0,0)$).

\subsection{R\'{e}duction du complexe
$\left(\!{\rm Hom}({\ve^\cdot\underline{{\gsI}_+^{\otimes \cdot}}},\gsI),
\left[[-,-]_S+\w,-\right]\right)$}

\noindent Dans cette partie on s'attache \`{a} d\'{e}montrer
\begin{prop}
Le complexe
$\left(\!{\rm Hom}({\ve^\cdot\underline{{\gsI}_+^{\otimes \cdot}}},\gsI),
\left[[-,-]_S\! +\! \w,-\right]\right)$ est acyclique
si le complexe $\left({\Hom}_{{\g}}(\ve_{{\g}}^\cdot
\g \ot_{{\gsI}_+}\Omega_{{\gsI}_+},\gsI),\delta^{1,1} \right)$ l'est.
\end{prop}
\begin{prooof}
La multiplication $m^2$ donne une structure de
$\g$-module \`{a} gauche \`{a} l'espace vectoriel $\g\ot V$
(par multiplication sur le premier facteur) pour tout espace gradu\'{e} V.
On a un isomorphisme de complexes~:
\begin{multline*}
\left({\rm Hom}(\gu,\gsI),[m^{1,1}+m^2,-]\right)\cong  \cr
\left({\rm Hom}_{{\g}}(\gextgsI,\gsI),[m^{1,1}+m^2,-]\right).
\end{multline*}
La diff\'{e}rentielle $[m^2,-]$ induite sur le dernier complexe est la duale
d'une diff\'{e}rentielle induite par
$\delta^2$ sur $\gextgsI$ qui n'est autre que la diff\'{e}rentielle de
Harrison $\beta$ sur chaque facteur
$\g\ot \ot \underline{\gsI^{\otimes \cdot}}$.
Un argument standard de suites spectrales~{\it c.f.}~\cite{Tam}, \cite{GH}
assure que si   l'homologie du complexe de Harrison
$\left(\g\ot \underline{{\gsI}_+^{\ot \cdot}}, \beta\right)$ de $\gsI$
\`{a} coefficient dans $\g$  est  \'{e}gale
\`{a} $\g \ot_{{\gsI}_+}\Omega^1_{{\gsI}_+}$ alors on a un quasi-isomorphisme
de complexes~:
$$\left({\Hom}_{{\g}}(\ve_{{\g}}^\cdot
\g \ot_{{\gsI}_+}\Omega_{{\gsI}_+},\gsI),\delta^{1,1} \right) \to
\left(\!{\rm Hom}({\ve^\cdot\underline{{\gsI}_+^{\otimes \cdot}}},\gsI),
\left[[-,-]_S+\w,-\right]\right)$$

On note  ${\gsI}_+=k\oplus \gsI$ l'alg\`{e}bre unitaire obtenue en ajoutant
une unit\'{e} \`{a} l'id\'{e}al $\gsI$.
Rappelons qu'en caract\'{e}ristique $0$, l'homologie de Harrison est \'{e}gale
\`{a} l'homologie d'Andr\'{e}-Quillen (\`{a} un d\'{e}calage du degr\'{e} de un pr\`{e}s)
que l'on note $AQ_\cdot(B/A,M)$ pour une $A$-alg\`{e}bre $B$ et un $A$-module
$M$.
On a une suite d'inclusions de sous-alg\`{e}bres gradu\'{e}es commutatives et
unitaires
$ k\hookrightarrow {\gsI}_+ \hookrightarrow \g $. Pour tout $\g$-module $M$,
la suite exacte de Jacobi-Zariski associ\'{e}e s'\'{e}crit
\begin{multline*}
\cdots AQ_{\cdot+1}(\g/{\gsI}_+, M)\to AQ_{\cdot}({\gsI}_+/k, M)\\
\to AQ_{\cdot}(\g/k, M)\to AQ_{\cdot}(\g/{\gsI}_+, M)\to \cdots .
\end{multline*}
On s'int\'{e}resse au cas $M=\g$.
On a un isomorphisme $AQ_{\cdot}(\g/k,\g)\cong \Omega^1_{\g/k}$ car $\g$ est
sym\'{e}trique ce qui ram\`{e}ne le calcul de $AQ_{\cdot}({\gsI}_+/k, \g)$
\`{a} celui de $AQ_{\cdot+1}(\g/{\gsI}_+, \g)$.
Rappelons que
$AQ_{\cdot}(\g/{\gsI}_+, M)\cong HH_{\cdot+1}^{(1)}(\g/{\gsI}_+, M)$
o\`{u} $HH_{\cdot}^{(n)}(B/A,M)$
d\'{e}signe la partie de poids $n$
dans la $\lambda$-d\'{e}composition de l'homologie de Hochschild
pour toute $A$-alg\`{e}bre $B$ et $B$-module $M$ (\cite{Lo}, Section~4).
Il est bien connu que l'homologie de Hochschild
est \'{e}gale \`{a} l'homologie du complexe normalis\'{e} d\'{e}fini,
pour tout $m\geq 0$, par
$\overline{C_m}^{{{\gsI}_+}}(\g, M)=M\ot_{{\gsI}_+}(\g/{{\gsI}_+})^{\ot_{{\gsI}_+}n}$
muni de la diff\'{e}rentielle de Hochschild $\mathbf{b}$.
On s'int\'{e}resse au cas $M=\g$.
Comme ${\gsI}_+=k\oplus \gsI$ et que $\gsI$ est un id\'{e}al,
la projection $\g\to \gsP$ induit un isomorphisme
\begin{align*}
\g\! \ot_{{\gsI}_+}\! \left(\g/{{\gsI}_+}\right)
\ot_{{\gsI}_+}\! \cdots \ot_{{\gsI}_+}\!
\left(\g/{{\gsI}_+}\right)&\cong \g/\gsI\ot_{k} \left(\gsP/k\right)
\ot_{k}\cdots \ot_{k} \left(\gsP/k\right)\cr
&\cong \gsP\ot_{k} \left(\gsP/k\right)\ot_{k}
\cdots \ot_{k} \left(\gsP/k\right).
\end{align*}
Le Th\'{e}or\`{e}me  d'Hochschild-Kostant-Rosenberg
appliqu\'{e} \`{a}
l'alg\`{e}bre sym\'{e}\-trique
$\gsP$ donne un isomorphisme
$$H_{\cdot}(\overline{C_\cdot}^{{k}}(\gsP,\gsP))\cong
HH_{\cdot}(\gsP)\cong \Lambda^{\cdot}\Omega^1_{\gsP}.$$
On en d\'{e}duit l'isomorphisme cherch\'{e}
$AQ_{\cdot}({\gsI}_+/k, \g)\cong \g\ot_{{\gsI}_+}\Omega^1_{{\gsI}_+}$ et
 la suite exacte de $\g$-modules
 $$\g\ot_{{\gsI}_+}\Omega^1_{{\gsI}_+}\hookrightarrow
 \Omega^1_{\g}\twoheadrightarrow \Omega^1_{\gsP}  $$
 o\`{u} la structure de $\g$-module de $\Omega^1_{\gsP}  $ est induite par la
 projection  $\g\twoheadrightarrow \gsP$.

\smallskip

\end{prooof}

\subsection{Acyclicit\'{e} du complexe $\left({\Hom}_{{\g}}(\ve_{{\g}}^\cdot
\g \ot_{{\gsI}_+}\Omega_{{\gsI}_+},\gsI),\delta^{1,1} \right)$}

Un calcul simple nous donne d\'{e}j\`{a}~:
\begin{prop}
Dans le cas $n=1=\nu$, le complexe
$\left({\Hom}_{{\g}}(\ve_{{\g}}^\cdot
\g \ot_{{\gsI}_+}\Omega_{{\gsI}_+}\right.$
$\left.,\gsI),\delta^{1,1} \right)$
est acycli\-que, c'est-\`{a}-dire qu'il est concentr\'{e} en degr\'{e} $0$.
\end{prop}

\medskip

Nous conjecturons que ce complexe est acyclique dans le cas
$\nu=1$: une indication forte est le fait qu'il est toujours possible de
repr\'{e}senter un star-produit, ce qui a
\'{e}t\'{e} montr\'{e} par Gl\"{o}{\ss}ner (voir \cite[Lemma 1]{Glo98}, \cite{Bor03}):
\begin{theorem}[Gl\"{o}{\ss}ner 1998]
Si $C$ est une sous-vari\'{e}t\'{e} co\"{\i}sotrope de codimension $1$ dans
$X$ et $*$ un star-produit sur $X$.
Alors on peut construire une star-repr\'{e}sentation.
\end{theorem}
L'acyclicit\'{e} du complexe ci-dessus nous permettrait,
une globalisation du th\'{e}or\`{e}me
de formalit\'{e} dans le cas $\nu=1$ (en reprenant la preuve de \cite{Dol}).

Enfin, ce travail nous donne une expression simple des obstructions
\`{a} la formalit\'{e} qui r\'{e}sident surtout dans la ou les diff\'{e}rentielles $d'$
possibles: on peut imaginer, aux vues des obstructions \`{a} la repr\'{e}sentabilit\'{e}
li\'{e}es aux classes
d'Atiyah-Molino jusqu'\`{a} l'ordre $3$ d'un star-produit symplectique
(voir \cite{Bor03}), qu'il
faille demander \`{a} la structure de Poisson compatible
des conditions additionnelles pour qu'elle soit repr\'{e}sentable.


\begin{thebibliography}{99}

\bibitem{BFFLS78}
  Bayen, F., Flato, M., Fr{\o}nsdal, C., Lichnerowicz, A., Sternheimer, D.:
  {\em Deformation Theory and Quantization.}
  Annals of Physics {\bf 111} (1978), part I: 61-110, part II: 111-151.

\bibitem{Bor03} Bordemann, M.: {\em (Bi)modules, morphismes et r\'{e}ductions
        des star-produits: le cas symplectique et classes caract\'{e}ristiques
        des feuilletages}. Pr\'{e}publication, \`{a} para\^{\i}tre.

\bibitem{BHW00}
     Bordemann, M., Herbig, H.-C., Waldmann, S.: {\em BRST cohomology and
     Phase Space
     Reduction in Deformation Quantisation}, Commun.Math.Phys. {\bf 210}
     (2000), 107-144.



\bibitem{Del} Deligne, P.: {\sl Letter to Stasheff},
Gerstenhaber May Schechtman, Drinfeld (1993).

\bibitem{Dol} Dolgushev, V.:
{\sl Covariant and Equivariant Formality Theorems}, I preprint
QA$\backslash$ 0307212 (2003).

\bibitem{EK} Etingof, P.,  Kazhdan, D.:
{\sl Quantization of Lie bialgebras I},
Selecta Math.,  N.S. (2) {\bf n.1} (1996), 1-41.
{\sl Quantization of Lie bialgebras II}, Selecta Math., N.S.(4)
{\bf n.2} (1998), 213-231, 233-269.

\bibitem{Ger63} Gerstenhaber, M.: {\em The Cohomology Structure of an
    Associative Ring}. Ann. Math. {\bf 78} (1963), 267-288.

\bibitem{GV} Gerstenhaber, M., Voronov, A.:
  {\it  Homotopy $G$-algebras and moduli space operad},
  Internat. Math. Res. Notices (1995), no. 3, 141--153

\bibitem{GH} Ginot, G., Halbout, G.:
{\sl A deformed version of Tamarkin's formality Theorem},
      pr\'{e}publication de l'IRMA (2002).

\bibitem{Glo98} Gl\"{o}{\ss}ner, P.: {\em Star-Product Reduction for Coisotropic
     Submanifolds of Codimension 1}. Pr\'{e}publication Facult\'{e} de Physique de
     l'Universit\'{e} de Freiburg FR-THEP-98/10, {\tt math.QA/9805049}, mai
     1998.
\bibitem{Kad02} Kadeishvili, T.: {\em Structure of $A(\infty)$-algebra
     and Hochschild and Harrison cohomology}. Proc.~of A.Razmadze
     Math.Inst. {\bf 91} (1988), 20-27, voir aussi {\tt math.AT/0210331}

\bibitem{Kon} Kontsevich, M.:
{\sl Deformation quantization of Poisson manifolds I}, pr\'{e}publication
IHES, {\tt QA/9709070} (1997).

\bibitem{Lo} Loday, J.-L.: {\sl Cyclic homology},
Springer-Verlag 1992.

\bibitem{Lu93} Lu, J.-H.: {\em Moment Maps at the Quantum Level}.
    Commun.~Math.~Phys.~{\bf 157} (1993), 389-404.


\bibitem{Sta} Stasheff, J.D.:
{\sl On the homology associativity 1 and II},
Transactions of the AMS 108 (1963), 275-292, 293-312.


\bibitem{Tam} Tamarkin, D.:
{\sl Another proof of M. Kontsevich formality theorem},
Preprint math$\backslash$9803025.


\end{thebibliography}
\end{document}